\documentclass[12pt]{article}
\usepackage[isolatin]{inputenc}
\usepackage{amsmath,amssymb,amsthm}
\usepackage{graphicx,pstricks}
\usepackage{subfigure}

\title{Spinning rough disk moving in a rarefied medium}
\author{Alexander Plakhov\thanks{Aberystwyth University, Aberystwyth SY23 3BZ, Ceredigion, UK and University of Aveiro, Aveiro 3810-193, Portugal} \and Tatiana Tchemisova \thanks{University of Aveiro, Aveiro 3810-193, Portugal} \and Paulo Gouveia \thanks{Polytechnic Institute of Bragança, Bragança 5301-854, Portugal}}

\date{}
      \newcommand {\al}   {\alpha}          \newcommand {\bt}  {\beta}
      \newcommand {\gam } {\gamma}          
      \newcommand {\del}  {\delta}          
              \newcommand {\ve}   {\varepsilon}
                 \newcommand {\vphi} {\varphi}
      \newcommand {\lam}  {\lambda}         
      
      \newcommand {\om}   {\omega}          \newcommand {\Om}  {\Omega}
      \newcommand {\pl}   {\partial}        
           
      \newcommand {\RRR}  {{\mathbb R}}

             \newcommand {\DDD}  {D}
      \newcommand {\rrr}  {\varrho}
     \newcommand {\beq}  {\begin{equation}}
      \newcommand {\eeq}  {\end{equation}}
     \newcommand {\beqa}  {\begin{equation*}}
      \newcommand {\eeqa}  {\end{equation*}}
    \newcommand {\interval}{[-\pi/2,\, \pi/2]}
      \newcommand{\Pigam} {\Pi}   \newcommand{\relmom}{\bt}
      \newcommand{\massgas}{M_g}
      
      \newtheorem{theorem}{Theorem}

\begin{document}
\maketitle

\begin{abstract}
We study the Magnus effect: deflection of the trajectory of a spinning body moving in a medium. The body is {\it rough}; that is, there are small cavities on its surface. We concentrate on the extreme case of rare medium, where  mutual interaction of the medium particles is neglected
and reflections of the particles from the body's surface are elastic. Now the effect is due solely to multiple reflections of particles in the cavities and can be studied by means of billiard theory.

We limit ourselves to the (two-dimensional) case of rough disk moving on the plane, and to a zero-temperature medium. We calculate the force of resistance acting on the body, as well as the moment of the force, and examine how these values depend on the ``kind of roughness'' (that is, on the shape of cavities). We show that a nonzero transversal component of the resistance force generally appears, resulting in deflection of the disk trajectory. The trajectory is determined for several particular simple cases.

We compare our results with the previous works on Magnus effect in highly rarefied gases, which derive the effect from {\it non-elastic} interaction of gas particles with the body. We propose another mechanism of creating the transversal force due to {\it multiple} collisions of particles with the body.
\\ \ \\

\end{abstract}

\begin{quote}
\end{quote}

\begin{quote}
\end{quote}

\begin{quote}
{\small {\bf Running title:} Spinning rough disk}
\end{quote}

\section{Introduction}

{\bf 1.1.} We are concerned with a spinning solid body moving in a homogeneous medium. The medium is extremely rarefied, so that the free path length of the medium particles is much larger than the body's size.
In such a case, the interaction of the body with the medium can be described in terms of {\it free molecular flow}, where a flow of point particles falls on the body's surface; each particle interacts with the body but does not with other particles. There is no gravitation force. We suppose that the particles of the medium stay at rest, that is, the absolute temperature of the medium equals zero.
In a frame of reference moving forward together with the body, we have a parallel flow of particles falling on the body  at rest.

Neglecting the angular momentum of particles, each particle is identified
with a mass point that approaches the body, makes several (maybe none) collisions with its surface, and then goes away. All reflections are supposed to be {\it perfectly elastic}.

In this paper, we restrict ourselves to the two-dimensional case. Consider a {\it rough circle}, that is, a set obtained from a circle by making infinitely small dimples on its boundary. More precisely, consider a sequence of sets $B_m$,\, $m = 3,\,4,\, 5, \ldots$, inscribed in the circle $B_r(O)$ of radius $r$ centered at  a point $O$. Each set $B_m$ is invariant under the rotation by the angle $2\pi/m$, and the intersection of $B_m$ with the $(2\pi/m)$-sector $A_m O  C_m$ formed by two radii $O A_m$ and $O  C_m$ is a set bounded by these radii and by a piecewise smooth non-self-intersecting curve contained in the triangle $A_m O  C_m$ and joining the points $A_m$ and $ C_m$. These curves are similar for all $m$. A {\it rough circle} (or {\it rough disk}) $B$, is associated with such a sequence  of sets $B_m$. The curve is called the {\it shape of roughness}. All the values related to  resistance or dynamics of the rough disk that are calculated below are understood as limits of the sequences of the corresponding values for $B_m$ as $m \to \infty$.

The roughness introduced here is uniform: it is identical at each point of the circle boundary.

The rough circle is rigid and can only undergo translations or rotations. Denote by $\vphi(t)$ the rotation angle at the time $t$, by $\om(t)$ the angular velocity of the body, $\om(t) = d\vphi/dt$, and let $\vec v(t)$ be the velocity of the body's center of mass. Let us agree to measure the rotation angle and the angular velocity counterclockwise.

\vspace{7mm}


\begin{figure}[h]
\begin{picture}(0,250)
   \rput(3,3.3){
   \scalebox{0.9}{
 \rput(0,-0.5){
 \psline[arrows=->,arrowscale=1.5](4,6.5)(4,5)
 \psline[arrows=->,arrowscale=1.5](3,6.5)(3,5)
 \psline[arrows=->,arrowscale=1.5](5,6.5)(5,5)
 \psline[arrows=->,arrowscale=1.5](6,6.5)(6,5)
 \psline[arrows=->,arrowscale=1.5](2,6.5)(2,5)
 \psline[arrows=->,arrowscale=1.5](1,6.5)(1,5)
 \psline[arrows=->,arrowscale=1.5](7,6.5)(7,5)
 }
 \rput(4,0){
\rput{-20}(-0.067,0.03){ \rput{-20}(-0.067,0.03){ \rput{-20}(-0.067,0.03){ \rput{-20}(-0.067,0.03){ \rput{-20}(-0.067,0.03){ \rput{-20}(-0.067,0.03){ \rput{-20}(-0.067,0.03){ \rput{-20}(-0.072,0.03){
\psecurve[linecolor=lightgray,linewidth=0pt,fillstyle=solid,fillcolor=lightgray](3.4,-0.1)(0,0)(3,1.27)(3.383,1.231)(3.2,1)(3,0.8)(3.4,0.4)(3.3,0.2)(3.4,0.1)(3.6,0)(3.4,-0.1)(0,0)(3.2,1.45)
\psecurve[linewidth=0.4pt](2.8,1.245)(3.06,1.285)(3.383,1.231)(3.2,1)(3,0.8)(3.4,0.4)(3.3,0.2)(3.32,0.165)(3.4,0.1)
} } } } } } } }
\rput{-20}(-0.067,0.03){ \rput{-20}(-0.067,0.03){ \rput{-20}(-0.067,0.03){ \rput{-20}(-0.067,0.03){
\rput{-20}(-0.067,0.03){ \rput{-20}(-0.067,0.03){ \rput{-20}(-0.072,0.03){
\psecurve[linecolor=lightgray,linewidth=0pt,fillstyle=solid,fillcolor=lightgray](3.4,-0.1)(0,0)(3,1.27)(3.383,1.231)(3.2,1)(3,0.8)(3.4,0.4)(3.3,0.2)(3.4,0.1)(3.6,0)(3.4,-0.1)(0,0)(3.2,1.45)
\psecurve[linewidth=0.4pt](2.8,1.245)(3.06,1.285)(3.383,1.231)(3.2,1)(3,0.8)(3.4,0.4)(3.3,0.2)(3.32,0.165)(3.4,0.1)
} } } } } } }
\rput{-20}(-0.067,0.03){ \rput{-20}(-0.067,0.03){ \rput{-20}(-0.067,0.03){ \rput{-20}(-0.067,0.03){ \rput{-20}(-0.067,0.03){ \rput{-20}(-0.072,0.03){
\psecurve[linecolor=lightgray,linewidth=0pt,fillstyle=solid,fillcolor=lightgray](3.4,-0.1)(0,0)(3,1.27)(3.383,1.231)(3.2,1)(3,0.8)(3.4,0.4)(3.3,0.2)(3.4,0.1)(3.6,0)(3.4,-0.1)(0,0)(3.2,1.45)
\psecurve[linewidth=0.4pt](2.8,1.245)(3.06,1.285)(3.383,1.231)(3.2,1)(3,0.8)(3.4,0.4)(3.3,0.2)(3.32,0.165)(3.4,0.1)
} } } } } }
\rput{-20}(-0.067,0.03){ \rput{-20}(-0.067,0.03){ \rput{-20}(-0.067,0.03){ \rput{-20}(-0.067,0.03){ \rput{-20}(-0.072,0.03){
\psecurve[linecolor=lightgray,linewidth=0pt,fillstyle=solid,fillcolor=lightgray](3.4,-0.1)(0,0)(3,1.27)(3.383,1.231)(3.2,1)(3,0.8)(3.4,0.4)(3.3,0.2)(3.4,0.1)(3.6,0)(3.4,-0.1)(0,0)(3.2,1.45)
\psecurve[linewidth=0.4pt](2.8,1.245)(3.06,1.285)(3.383,1.231)(3.2,1)(3,0.8)(3.4,0.4)(3.3,0.2)(3.32,0.165)(3.4,0.1)
} } } } }
\rput{-20}(-0.067,0.03){ \rput{-20}(-0.067,0.03){ \rput{-20}(-0.067,0.03){ \rput{-20}(-0.072,0.03){
\psecurve[linecolor=lightgray,linewidth=0pt,fillstyle=solid,fillcolor=lightgray](3.4,-0.1)(0,0)(3,1.27)(3.383,1.231)(3.2,1)(3,0.8)(3.4,0.4)(3.3,0.2)(3.4,0.1)(3.6,0)(3.4,-0.1)(0,0)(3.2,1.45)
\psecurve[linewidth=0.4pt](2.8,1.245)(3.06,1.285)(3.383,1.231)(3.2,1)(3,0.8)(3.4,0.4)(3.3,0.2)(3.32,0.165)(3.4,0.1)
} } } }
\rput{-20}(-0.067,0.03){
\rput{-20}(-0.067,0.03){
\rput{-20}(-0.072,0.03){
\psecurve[linecolor=lightgray,linewidth=0pt,fillstyle=solid,fillcolor=lightgray](3.4,-0.1)(0,0)(3,1.27)(3.383,1.231)(3.2,1)(3,0.8)(3.4,0.4)(3.3,0.2)(3.4,0.1)(3.6,0)(3.4,-0.1)(0,0)(3.2,1.45)
\psecurve[linewidth=0.4pt](2.8,1.245)(3.06,1.285)(3.383,1.231)(3.2,1)(3,0.8)(3.4,0.4)(3.3,0.2)(3.32,0.165)(3.4,0.1)
} } }
\rput{-20}(-0.067,0.03){
\rput{-20}(-0.072,0.03){
\psecurve[linecolor=lightgray,linewidth=0pt,fillstyle=solid,fillcolor=lightgray](3.4,-0.1)(0,0)(3,1.27)(3.383,1.231)(3.2,1)(3,0.8)(3.4,0.4)(3.3,0.2)(3.4,0.1)(3.6,0)(3.4,-0.1)(0,0)(3.2,1.45)
\psecurve[linewidth=0.4pt](2.8,1.245)(3.06,1.285)(3.383,1.231)(3.2,1)(3,0.8)(3.4,0.4)(3.3,0.2)(3.32,0.165)(3.4,0.1)
} }
\rput{-20}(-0.072,0.03){
\psecurve[linecolor=lightgray,linewidth=0pt,fillstyle=solid,fillcolor=lightgray](3.4,-0.1)(0,0)(3,1.27)(3.383,1.231)(3.2,1)(3,0.8)(3.4,0.4)(3.3,0.2)(3.4,0.1)(3.6,0)(3.4,-0.1)(0,0)(3.2,1.45)
\psecurve[linewidth=0.4pt](2.8,1.245)(3.06,1.285)(3.383,1.231)(3.2,1)(3,0.8)(3.4,0.4)(3.3,0.2)(3.32,0.165)(3.4,0.1)
}
\rput{20}(-0.067,-0.03){ \rput{20}(-0.067,-0.03){ \rput{20}(-0.067,-0.03){ \rput{20}(-0.067,-0.03){
\rput{20}(-0.067,-0.03){ \rput{20}(-0.067,-0.03){ \rput{20}(-0.067,-0.03){ \rput{20}(-0.067,-0.03){
\rput{20}(-0.075,-0.03){
\psecurve[linecolor=lightgray,linewidth=0pt,fillstyle=solid,fillcolor=lightgray](3.4,-0.1)(0,0)(3,1.27)(3.383,1.231)(3.2,1)(3,0.8)(3.4,0.4)(3.3,0.2)(3.4,0.1)(3.6,0)(3.4,-0.1)(0,0)(3.2,1.45)
\psecurve[linewidth=0.4pt](2.8,1.245)(3.06,1.285)(3.383,1.231)(3.2,1)(3,0.8)(3.4,0.4)(3.3,0.2)(3.32,0.165)(3.4,0.1)
} } } } } } } } }
\rput{20}(-0.067,-0.03){ \rput{20}(-0.067,-0.03){ \rput{20}(-0.067,-0.03){ \rput{20}(-0.067,-0.03){
\rput{20}(-0.067,-0.03){ \rput{20}(-0.067,-0.03){ \rput{20}(-0.067,-0.03){ \rput{20}(-0.075,-0.03){
\psecurve[linecolor=lightgray,linewidth=0pt,fillstyle=solid,fillcolor=lightgray](3.4,-0.1)(0,0)(3,1.27)(3.383,1.231)(3.2,1)(3,0.8)(3.4,0.4)(3.3,0.2)(3.4,0.1)(3.6,0)(3.4,-0.1)(0,0)(3.2,1.45)
\psecurve[linewidth=0.4pt](2.8,1.245)(3.06,1.285)(3.383,1.231)(3.2,1)(3,0.8)(3.4,0.4)(3.3,0.2)(3.32,0.165)(3.4,0.1)
} } } } } } } }
\rput{20}(-0.067,-0.03){ \rput{20}(-0.067,-0.03){ \rput{20}(-0.067,-0.03){ \rput{20}(-0.067,-0.03){
\rput{20}(-0.067,-0.03){ \rput{20}(-0.067,-0.03){ \rput{20}(-0.075,-0.03){
\psecurve[linecolor=lightgray,linewidth=0pt,fillstyle=solid,fillcolor=lightgray](3.4,-0.1)(0,0)(3,1.27)(3.383,1.231)(3.2,1)(3,0.8)(3.4,0.4)(3.3,0.2)(3.4,0.1)(3.6,0)(3.4,-0.1)(0,0)(3.2,1.45)
\psecurve[linewidth=0.4pt](2.8,1.245)(3.06,1.285)(3.383,1.231)(3.2,1)(3,0.8)(3.4,0.4)(3.3,0.2)(3.32,0.165)(3.4,0.1)
} } } } } }}
\rput{20}(-0.067,-0.03){ \rput{20}(-0.067,-0.03){ \rput{20}(-0.067,-0.03){ \rput{20}(-0.067,-0.03){
\rput{20}(-0.067,-0.03){ \rput{20}(-0.075,-0.03){ \psecurve[linecolor=lightgray,linewidth=0pt,fillstyle=solid,fillcolor=lightgray](3.4,-0.1)(0,0)(3,1.27)(3.383,1.231)(3.2,1)(3,0.8)(3.4,0.4)(3.3,0.2)(3.4,0.1)(3.6,0)(3.4,-0.1)(0,0)(3.2,1.45)
\psecurve[linewidth=0.4pt](2.8,1.245)(3.06,1.285)(3.383,1.231)(3.2,1)(3,0.8)(3.4,0.4)(3.3,0.2)(3.32,0.165)(3.4,0.1)
} } } } } }
\rput{20}(-0.067,-0.03){ \rput{20}(-0.067,-0.03){ \rput{20}(-0.067,-0.03){ \rput{20}(-0.067,-0.03){
\rput{20}(-0.075,-0.03){
\psecurve[linecolor=lightgray,linewidth=0pt,fillstyle=solid,fillcolor=lightgray](3.4,-0.1)(0,0)(3,1.27)(3.383,1.231)(3.2,1)(3,0.8)(3.4,0.4)(3.3,0.2)(3.4,0.1)(3.6,0)(3.4,-0.1)(0,0)(3.2,1.45)
\psecurve[linewidth=0.4pt](2.8,1.245)(3.06,1.285)(3.383,1.231)(3.2,1)(3,0.8)(3.4,0.4)(3.3,0.2)(3.32,0.165)(3.4,0.1)
} } } } }
\rput{20}(-0.067,-0.03){ \rput{20}(-0.067,-0.03){ \rput{20}(-0.067,-0.03){ \rput{20}(-0.075,-0.03){
\psecurve[linecolor=lightgray,linewidth=0pt,fillstyle=solid,fillcolor=lightgray](3.4,-0.1)(0,0)(3,1.27)(3.383,1.231)(3.2,1)(3,0.8)(3.4,0.4)(3.3,0.2)(3.4,0.1)(3.6,0)(3.4,-0.1)(0,0)(3.2,1.45)
\psecurve[linewidth=0.4pt](2.8,1.245)(3.06,1.285)(3.383,1.231)(3.2,1)(3,0.8)(3.4,0.4)(3.3,0.2)(3.32,0.165)(3.4,0.1)
} } } }
\rput{20}(-0.067,-0.03){ \rput{20}(-0.067,-0.03){ \rput{20}(-0.075,-0.03){
\psecurve[linecolor=lightgray,linewidth=0pt,fillstyle=solid,fillcolor=lightgray](3.4,-0.1)(0,0)(3,1.27)(3.383,1.231)(3.2,1)(3,0.8)(3.4,0.4)(3.3,0.2)(3.4,0.1)(3.6,0)(3.4,-0.1)(0,0)(3.2,1.45)
\psecurve[linewidth=0.4pt](2.8,1.245)(3.06,1.285)(3.383,1.231)(3.2,1)(3,0.8)(3.4,0.4)(3.3,0.2)(3.32,0.165)(3.4,0.1)
} } }
\rput{20}(-0.067,-0.03){
\rput{20}(-0.075,-0.03){
\psecurve[linecolor=lightgray,linewidth=0pt,fillstyle=solid,fillcolor=lightgray](3.4,-0.1)(0,0)(3,1.27)(3.383,1.231)(3.2,1)(3,0.8)(3.4,0.4)(3.3,0.2)(3.4,0.1)(3.6,0)(3.4,-0.1)(0,0)(3.2,1.45)
\psecurve[linewidth=0.4pt](2.8,1.245)(3.06,1.285)(3.383,1.231)(3.2,1)(3,0.8)(3.4,0.4)(3.3,0.2)(3.32,0.165)(3.4,0.1)
} }
\rput{20}(-0.075,-0.03){
\psecurve[linecolor=lightgray,linewidth=0pt,fillstyle=solid,fillcolor=lightgray](3.4,-0.1)(0,0)(3,1.27)(3.383,1.231)(3.2,1)(3,0.8)(3.4,0.4)(3.3,0.2)(3.4,0.1)(3.6,0)(3.4,-0.1)(0,0)(3.2,1.45)
\psecurve[linewidth=0.4pt](2.8,1.245)(3.06,1.285)(3.383,1.231)(3.2,1)(3,0.8)(3.4,0.4)(3.3,0.2)(3.32,0.165)(3.4,0.1)
}
\psecurve[linecolor=lightgray,linewidth=0pt,fillstyle=solid,fillcolor=lightgray](3.4,-0.1)(0,0)(3,1.27)(3.383,1.231)(3.2,1)(3,0.8)(3.4,0.4)(3.3,0.2)(3.4,0.1)(3.6,0)(3.4,-0.1)(0,0)(3.2,1.45)
\psecurve[linewidth=0.4pt](2.8,1.245)(3.06,1.285)(3.383,1.231)(3.2,1)(3,0.8)(3.4,0.4)(3.3,0.2)(3.32,0.165)(3.4,0.1)
}
 \pscircle[linestyle=dashed,linewidth=0.6pt](4.1,0){3.6}
 \psarc[linewidth=1pt,arrows=->,arrowscale=2](4,0){2.3}{50}{130}
 \rput(4,2){\Large $\om$}
 \rput(4,0){\scalebox{1.2}{\Huge $B$}}
  }}
\end{picture}
\label{fig rough circle.} \caption{A rotating rough disk in a parallel flow of particles.}
\end{figure}
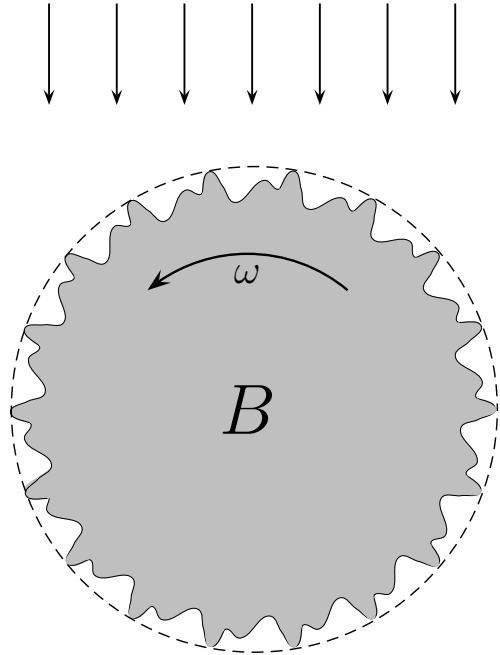
 \vspace{5mm}

Here we consider two problems as follows: (A) determine the force
of the medium resistance acting on the body, find the moment of this
force with respect to the body's center of mass, and investigate their
dependence on the shape of roughness, and (B) analyze the motion
of the body in the medium, that is, study the functions $\om(t)$ and $\vec v(t)$.
Problem (A) is primary with respect to problem (B). In the paper, we will
devote the main attention to problem (A), having just touched upon
problem (B); we will restrict ourselves to deducing equations of motion
and solving these equations for some simple particular cases.

With each set $B_m$ we associate a  distribution of mass inside $B_m$ such that the total mass $M$ is constant and the center of mass coincides with the center $O$ of the set. We also assume that the moment of inertia $I_m$ of $B_m$ about $O$ converges to a positive value $I$ as $m \to \infty$. One always has $I \le Mr^2$.  Denote $\relmom = Mr^2/I$;  that is, $\relmom$ is the inverse relative moment of inertia; we  have $1 \le \relmom < +\infty$. In what follows we will pay special attention to the two particular cases: (a) $\relmom = 1$; that is, the mass of the rough circle is concentrated near its boundary,\, and (b) $\relmom = 2$; the mass is distributed uniformly in the circle.

The resistance force $\vec R_m(B_m, \vphi, \om, \vec v)$ acting on the set $B_m$ and the moment of this force  $R_{I,m}(B_m, \vphi, \om, \vec v)$ depend on the shape of roughness, the rotation angle $\vphi$, the angular velocity $\om$, and the velocity of translation $\vec v$. The equations of dynamics are
 \begin{equation*}\label{dynamic 1}
M\, \frac{d\vec v}{dt}\, =\, \vec R_m(B_m, \vphi, \om, \vec v),
 \end{equation*}
 \begin{equation*}\label{dynamic 2}
I_m\, \frac{d\om}{dt}\, =\, R_{I,m}(B_m, \om, \vec v),
 \end{equation*}
 \begin{equation*}\label{dynamic 3}
\frac{d\vphi}{dt}\, =\, \om. \quad \quad \quad \quad \quad
 \end{equation*}
Taking the limit $m \to \infty$, one gets the resistance force $\vec R(B, \om, \vec v)$ and the moment of the force $R_I(B, \om, \vec v)$ of the rough disk, which do not depend on $\vphi$ anymore, and the equations for the disk dynamics take the form
  \beq\label{dynam1}
M\, \frac{d\vec v}{dt}\, =\, \vec R(B, \om, \vec v),
 \eeq
 \beq\label{dynam2}
I\, \frac{d\om}{dt}\, =\, R_I(B, \om, \vec v).
 \eeq \\

{\bf 1.2.} We shall see below that, generally speaking, $\vec R = \vec R(B, \om, \vec v)$
is not collinear to $\vec v$.

If a transversal force acts on a spinning body moving in a medium, resulting in deflection of the body's trajectory, then we encounter the so called {\it Magnus effect}.  Well-known examples of this phenomenon are dynamics of a golf ball and a football. The most intensively studied shapes are sphere and cylinder (in the second case, the cylinder is supposed to move in a direction orthogonal to the symmetry axis). If the direction of the transversal force coincides with the instantaneous velocity of the front point of the body, then a proper Magnus effect takes place. If these directions are opposite, then a {\it reverse} Magnus effect occurs. See Fig.\,2\,(a),\,(b).

\vspace{10mm}
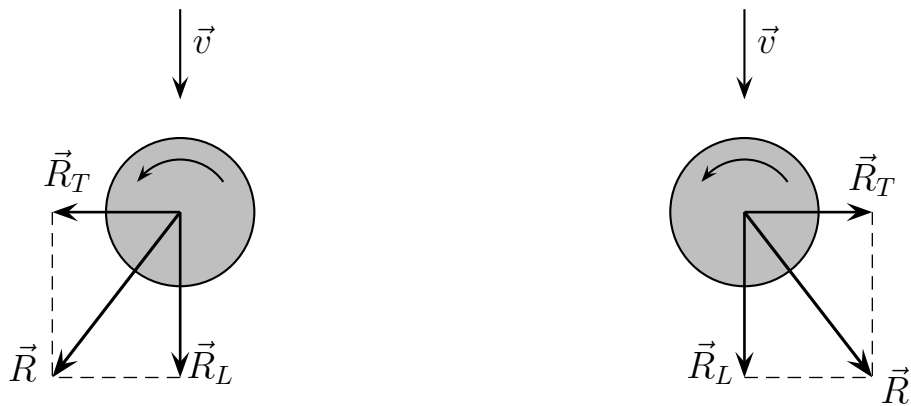
\begin{figure}[h]
\begin{picture}(0,160)
  \rput(0.5,2.7){
 \rput(10,0){
 \pscircle[fillstyle=solid,fillcolor=lightgray](0,0){1}
 \psline[arrows=<-,arrowscale=2](0,1.5)(0,2.7)
 \rput(0.3,2.3){\large $\vec v$}
 \psarc[linewidth=0.8pt,arrows=->,arrowscale=1.5](0,0){0.7}{35}{145}
 \psline[linewidth=1pt,arrows=->,arrowscale=2](0,0)(0,-2.2)
 \psline[linewidth=1pt,arrows=->,arrowscale=2](0,0)(1.7,0)
 \psline[linewidth=1.2pt,arrows=->,arrowscale=2](0,0)(1.7,-2.2)
 \psline[linewidth=0.5pt,linestyle=dashed](0,-2.2)(1.7,-2.2)(1.7,0)
 \rput(2,-2.3){\large $\vec R$}
 \rput(1.7,0.5){\large $\vec R_T$}
 \rput(-0.45,-2){\large $\vec R_L$}
 }
 \rput(2.5,0){
 \pscircle[fillstyle=solid,fillcolor=lightgray](0,0){1}
 \psline[arrows=<-,arrowscale=2](0,1.5)(0,2.7)
 \rput(0.3,2.3){\large $\vec v$}
 \psarc[linewidth=0.8pt,arrows=->,arrowscale=1.5](0,0){0.7}{35}{145}
 \psline[linewidth=1pt,arrows=->,arrowscale=2](0,0)(0,-2.2)
 \psline[linewidth=1pt,arrows=->,arrowscale=2](0,0)(-1.7,0)
 \psline[linewidth=1.2pt,arrows=->,arrowscale=2](0,0)(-1.7,-2.2)
 \psline[linewidth=0.5pt,linestyle=dashed](0,-2.2)(-1.7,-2.2)(-1.7,0)
 \rput(-2.1,-2){\large $\vec R$}
 \rput(-1.5,0.55){\large $\vec R_T$}
 \rput(0.4,-2){\large $\vec R_L$}
  }}
\end{picture}
\label{fig Magnus effect} \caption{(a) Magnus effect. \hspace*{20mm} (b) Reverse Magnus effect.}
\end{figure}
\vspace{7mm}

 There is a vast literature devoted to Magnus effect, motivated by sports and technology applications (see, e.g., \cite{Prandtl},\cite{RK},\cite{Mehta}). Under normal conditions, the {\it direct} Magnus effect was found to take place. However, in certain circumstances (for example, a well polished ball) the {\it reverse} Magnus effect can occur.

Nowadays, the growing interest to aerodynamics of artificial satellites stimulates the study of this phenomenon in {\it rarefied} media \cite{BSE},\cite{IY},\cite{W},\cite{WH}. In these papers the medium is supposed to be so rarefied that a description in terms of free molecular flow is possible. They suggest that in this case the {\it reverse} Magnus effect occurs, and study this phenomenon. Moreover, in a recent paper the effect of Magnus force on a spinning artificial satellite is investigated \cite{BS}.

In our opinion, the reverse Magnus effect in  highly rarefied media is caused by two factors:

(i) Non-elastic interaction of particles with the body. A part of the tangential component of the momentum of particles is transmitted to the body, resulting in creation of a transversal force.

(ii) Multiple collisions of particles with the body due to the fact that the body's surface is not convex but contains microscopic cavities.

In the papers \cite{BSE},\cite{IY},\cite{W},\cite{WH} the impact of factor (i) is studied. Moreover, in these works the body is supposed to be convex and therefore factor (ii) is excluded. In the present paper, in contrast, we concentrate on the study of factor (ii). We suppose that all collisions are perfectly elastic and so exclude the influence of factor (i).

We believe that the methods developed here can be extended to the three-dimensional case, and to the case of a finite temperature medium. Later on we shall see that although both reverse and direct Magnus effects are possible in our model, the reverse effect is the predominant one.

Note that the limiting case of slow rotation was considered in \cite{MatSb 04: mean resist},\cite{Arch Rat Mech},\cite{rough body 2D}. Actually, the class of bodies studied there is rather large and includes both rough disks and bounded simply connected sets. It was found that the mean resistance force is parallel to the direction of the body's motion, and therefore the Magnus effect does not appear.

In the next section, to each rough circle $B$ we assign a measure $\nu_B$ characterizing the law of billiard scattering on $B$. The values $\vec R$ and $R_I$ are defined to be functions of the values $\nu_B, \ \om$, and $\vec v$. In section 3 we calculate $\vec R$ and $R_I$ for some special values of $\nu_B$ (and thus for some special kinds of rough bodies). Further, we define the set of all possible values of $\vec R$, when $\om$ and $\vec v$ are fixed and $\nu_B$ takes all  admissible values. In other words, we answer the following question: what is the range of values for the resistance force acting on a rough unit circle? While we look through all possible shapes of roughness, the values of vector $\vec R$ cover a fixed convex two-dimensional set. The problem of determination of this set is formulated in terms of a special vector-valued Monge-Kantorovich problem and is solved numerically for some fixed values of the parameter $\lam = \om r/v$. In section 4 we deduce the equations of dynamics in a convenient form and solve these equations for some simple particular cases. Finally, in section 5 a comparison of our results with the previous works on reverse Magnus effect in rarefied media is given.

\section{Law of billiard scattering and resistance}

\subsection{Billiard scattering on a nonconvex set}

Let us first define the measure $\nu_\DDD$ characterizing billiard scattering on a bounded simply connected set $\DDD   \subset \RRR^2$ with piecewise smooth boundary.
The set $\pl (\text{conv} \DDD) \setminus \pl \DDD = \cup_{i \ge 1} I_i$ is the union of a finite or countable family of connected components $I_i, \ i = 1,\, 2, \ldots$. Each component $I_i$ is an open interval. (On Fig.\,3, the intervals $I_i$,\, $i = 1,\, 2,\, 3$ are shown dashed.) Denote $I_0 = \pl(\text{conv} \DDD) \cap \pl \DDD$; in other words, $I_0$ is the ``convex part'' of the boundary $\pl \DDD$. Thus, $\pl \DDD$ is
the disjoint union $\pl \DDD = \cup_{i \ge 0} I_i$.
 \vspace{10mm}

\begin{figure}[h]
\begin{picture}(0,185)
 \rput(1.7,3){
\scalebox{0.9}{
 \rput(5.8,0){\psecurve[fillstyle=solid,fillcolor=lightgray](3,0.5)(3.4,0)(3,-0.5)
 (3.3,-1.1)(2.5,-1.7) (2.3,-2.6)(1.9,-2.4)(1.3,-3.3)(0.4,-3.1)
 (0.3,-3.4)(-0.5,-3)(-1.3,-3.3)(-1.8,-2.7)(-2.45,-2.45)(-2.5,-1.7)(-3.15,-1.3)(-3.45,-0.8)(-3.35,0.2)(-3.1,0.2)
 (-3.1,1.6)(-2.8,1.6)(-2.05,2.5)(-1.7,2.6)(-1.3,3.3)(-0.6,3.1)(0.3,3.5)(0.2,3.2)(1.5,3.2)(1.8,2.6)
 (2.1,2.7)(3.0,1.8)(3.3,1.3)(3.1,0.6)(3.4,0)(3,0.5)
 \psline[linestyle=dashed](-3.51,-0.05)(-3.18,1.47)
 \rput(-3.55,0.7){$I_1$}
  \rput(-2.75,0.7){$\Om_1$}
 \psline[linestyle=dashed](-3.12,1.6)(-1.4,3.23)
 \rput(-2.5,2.55){$I_2$}
  \rput(-1.85,2.25){$\Om_2$}
 \psline[linestyle=dashed](-1.25,3.32)(0.15,3.55)
 \rput(-0.6,3.75){$I_3$}
  \rput(-0.5,2.85){$\Om_3$}
     \psdots(-3.08,-1.5)(-2.65,-2.15)   
     \psline[linewidth=0.5pt,arrows=->,arrowscale=1.5](-4.2,-1.74)(-3.53,-1.595)
    \psline[linewidth=0.5pt,arrows=->,arrowscale=1.5](-3.53,-1.595)(-2.86,-1.45)(-2.48,-2.0)(-3.5,-3)
  \psline[linestyle=dashed,linewidth=0.5pt](-3.22,-1.27)(-2.53,-2.33)
    \rput(-3.4,-1.3){$\xi$}
    \rput(-2.7,-2.6){$\xi^+$}
    \psline[linestyle=dotted](-3.08,-1.5)(-4.14,-2.19)
     \psline[linestyle=dotted](-2.65,-2.15)(-3.71,-2.84)
     \psarc(-3.08,-1.5){0.7}{192.21}{213.06}
     \psarc(-2.65,-2.15){0.8}{213.06}{224.43}
     \psarc(-2.65,-2.15){0.73}{213.06}{224.43}
       \rput(-2.15,-1.8){$\Om_i$}
    \rput(-4,-1.9){$\vphi$}
   \rput(-3.3,-2.4){$\vphi^+$}
    \rput(0,0){\Huge $D$}
     }
}}
\end{picture}
\label{fig cavities} \caption{Cavities on the boundary of a nonconvex set.}
\end{figure}
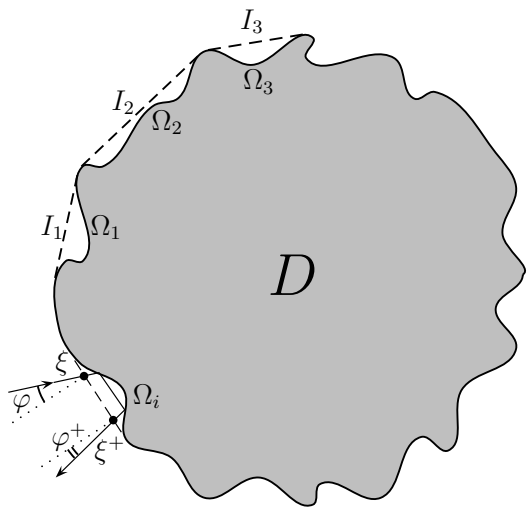

Further, the set conv$\DDD \setminus \DDD$ is the union of a finite or countable collection of its connected components. For any $I_i$ there exists a set $\Om_i$ from this collection such that $I_i \subset \pl\Om_i$; see Fig.\,3. The pair $(\Om_i, I_i)$ will be called a {\it cavity}, and the interval $I_i$, {\it
the opening of the cavity}.

Denote by $\vec n_\xi$ the outer unit normal to $\pl (\text{conv} \DDD)$ at the point $\xi \in \pl (\text{conv} \DDD)$. On the set $I_i \times \interval$ with the coordinates $(\xi, \vphi)$, define the measure $\mu_i$ according to the formula
$d\mu_i = \cos\vphi\, d\xi\, d\vphi$, where $d\xi$ and $d\vphi$ stand for the one-dimensional Lebesgue measure. Consider the billiard in $\RRR^2 \setminus \DDD$. Fix $\xi \in I_i$ and $\vphi \in \interval$, and consider the billiard particle that starts to move from the point $\xi$ with the velocity forming the angle $\vphi$ with $-\vec n_\xi$. The particle makes one or several reflections at points of $\pl\Om_i \setminus I_i$ and then  intersects $I_i$ once {again} at the point $\xi^+ = \xi^+_i(\xi, \vphi)$, the velocity at the moment of intersection forming the angle $\vphi^+ = \vphi^+_i(\xi, \vphi)$ with the vector $\vec n_{\xi^+},\,$ $ \ \vphi^+ \in \interval$. Notice that it is always valid $\vec n_{\xi} = \vec n_{\xi^+}$. On Fig.\,3, we have $\vphi < 0$ and $\vphi^+ = \vphi_i^+(\xi, \vphi) > 0$.

Thus, for each $i$ we have defined the mapping $(\xi, \vphi) \mapsto (\xi_i^+, \vphi_i^+)$ from a full measure subset of $I_i \times \interval$ onto itself. This mapping is a bijection and involution, and preserves the measure $\mu_i$. In the particular case, where $i = 0$, it holds $\xi_0^+ = \xi$ and $\vphi_0^+ = -\vphi$. Denote by $l_i = |I_i|$ the length of $I_i$, and by $l = \sum_{i \ge 0} l_i = |\pl(\text{conv} \DDD)|$ the perimeter of $\text{conv} \DDD$. Introduce the notation $\Box :=  \interval \times $ $\interval$ and define the measures $\nu_\DDD^i$,\, $i = 0,\, 1,\, 2,\ldots$ on $\Box$ as follows: $\nu_\DDD^i(A) := \frac{1}{l_i} \mu_i \big( \left\{ (\xi, \vphi) : \right.$ $\left. (\vphi, \vphi_i^+(\xi, \vphi)) \in A \right\} \big)$ for any Borel set $A \subset \Box$. In particular, the measure $\nu_\DDD^0 =: \nu_0$ is supported on the diagonal $\vphi^+ = -\vphi$, and its density equals $d\nu_0(\vphi,\vphi^+) = \cos\vphi \cdot \delta(\vphi + \vphi^+)$. Finally, define $\nu_\DDD :=  \frac{1}{l} \sum_{i \ge 0} l_i \nu_\DDD^i$.

In a less formal way, the measure $\nu_\DDD$ can be interpreted as follows. Place the body $\DDD$ in a kind of ``ether'' and keep it motionless. Here by ether we mean a homogeneous isotropic medium composed of (mutually noninteracting) particles moving freely with unit velocity in all possible directions. When colliding with the body, the particles reflect elastically from its boundary. Thus, the particles of the ether behave like billiard particles in $\RRR^2 \setminus \DDD$. Homogeneity and isotropy mean that the ether density is constant and the velocities of the particles are uniformly distributed in $S^1$.

For each particle that has reflected from $\DDD$, fix the pair $(\vphi, \vphi^+)$ consisting of the angle of incidence $\vphi$ and the angle of reflection $\vphi^+$. The angle $\vphi$ is formed by the initial velocity of the particle and the vector $-\vec n_\xi$, whereas the angle $\vphi^+$ is formed by its final velocity and the vector $\vec n_{\xi^+}$. Here $\xi$ and $\xi^+$ are the points of the first and second intersection of the particle's trajectory with $\pl(\text{conv} \DDD)$. {\em The distribution of the set of pairs $(\vphi, \vphi^+)$ for all particles that have collided with $\DDD$ during a unit time period is described by the measure $\nu_\DDD$.}

Recall that given a Borel mapping $\pi: X \to Y$ of two sets $X \subset \RRR^{d_1}$,\, $Y \subset \RRR^{d_2}$, with the set $X$ being equipped with a Borel measure $\mu$, the so-called push-forward measure $\pi^\# \mu$ on $Y$ is defined  by $\pi^\# \mu(A) := \mu(\pi^{-1}(A))$ for any Borel set $A \subset Y$ . Denote by $\pi_{\vphi}$ and $\pi_{\vphi^+} : \Box \to \interval$ the projections of the square $\Box$ on its horizontal and vertical sides, respectively; $\pi_{\vphi}(\vphi, \vphi^+) = \vphi$,\, $\pi_{\vphi^+}(\vphi, \vphi^+) = \vphi^+$. The push-forward measures $\pi_\vphi^\# \nu$ and $\pi_{\vphi^+}^\# \nu$ are called marginal measures for $\nu$. They are defined on $\interval$, and for any Borel set $A \subset \interval$  holds
$$
\pi_{\vphi}^\# \nu (A) = \nu(A \times
\interval),\, \ \ \ \ \pi_{\vphi^+}^\# \nu (A) = \nu(\interval
\times A).
$$

Define the transformation of the square $\pi_d : \Box \to \Box$ exchanging the coordinates $\vphi$ and $\vphi^+$; that is, $\pi_d(\vphi, \vphi^+) = (\vphi^+, \vphi)$. The push-forward measure $\pi_d^\# \nu$ satisfies the condition $\pi_d^\# \nu(A) = \nu(\pi_d(A))$ for any Borel set $A \subset \Box$. In other words, the measures $\pi_d^\# \nu$ and $\nu$ are mutually symmetric with respect to the diagonal $\vphi = \vphi^+$.

Finally, define the measure $\gam$ on $\interval$ by $d\gam = \cos\vphi\, d\vphi$, and denote by $\Pigam^{\text{symm}}$ the set of measures $\nu$ on $\Box$ satisfying the conditions
 \beq\label{uslovie 1}
\pi_{\vphi}^\# \nu = \gam = \pi_{\vphi^+}^\# \nu
 \eeq
and
 \beq\label{uslovie 2}
\pi_d^\# \nu = \nu.
 \eeq
In other words, a generic measure from $\Pigam^{\text{symm}}$ is symmetric with respect to the diagonal $\vphi = \vphi^+$, and both of its marginal measures coincide with $\gam$. One has $\nu_\DDD^i \in \Pigam^{\text{symm}}$; this can be easily deduced from the measure preserving and involutive properties of the mapping $(\xi, \vphi) \mapsto (\xi_i^+, \vphi_i^+)$; for details see \cite{MatSb 04: mean resist}.  Hence $\nu_\DDD \in \Pigam^{\text{symm}}$. The following reverse statement also holds true.

\begin{theorem}\label{t1}
Whatever the two sets $K_1 \subset K_2 \subset \RRR^2$ such that {\rm dist}$(\pl K_1, \pl K_2) > 0$, the set of measures $\{ \nu_\DDD : \ K_1 \subset \DDD \subset K_2 \}$ is everywhere dense in $\Pigam^{\text{\rm symm}}$
in the weak topology.
\end{theorem}

The proof of this theorem can be found in \cite{Arch Rat Mech}.

\vspace{5 mm}

Now consider the rough circle $B$ generated by a sequence of sets $B_m$. All the cavities of all the sets $B_m$ are similar, therefore $\nu_{B_m}$ does not depend on $m$ and we can set by definition $\nu_B := \nu_{B_m}$. We will see later that  the resistance of $B$ can be written down as a functional of $\nu_B$. The following theorem is obtained by a slight modification of the proof of Theorem \ref{t1}.

\begin{theorem}\label{teor2}
Whatever $r > 0$, the set $\{ \nu_B : \ B \text{ is a rough circle of radius r} \, \}$ is everywhere dense in $\Pigam^{\text{\rm symm}}$ in the weak topology.
\end{theorem}

\subsection{Resistance of a rough disk}

Denote $v = |\vec v|$ and choose the  (non-inertial) frame of reference $Ox_1 x_2$ in such a way that the direction of the axis $Ox_2$ coincides with the direction of the  disk motion and the origin $O$ coincides with the  disk center.  In this frame of reference the  disk stays at rest, and the flow of particles falls down on it at the velocity $-\vec v_0 = (0;\, -v)^T$. Here and in what follows, we represent vectors as columns; for instance, a vector $\vec x$ will be denoted by $\left[ \!\! \begin{array}{c} x_1\\ x_2 \end{array} \!\! \right]$ or $(x_1;\, x_2)^T$.

 Let us calculate the force $\vec R$ of the medium resistance and the moment of this force $R_I$ with respect to $O$. To that end, first we consider the prelimit body $B_m$. Parameterize the opening of each cavity by the variable $\xi$ varying from 0 to 1 (recall that all the cavities are identical). Denote by $\rho$ the flow density, by $\vphi$, the rotation angle of the cavity (that is, the external normal at the cavity opening equals $\vec n_\xi = (-\sin\vphi;\, \cos\vphi)^T$), and by $\vec v^+_{(m)}(\xi, \vphi)$, the final velocity of the particle entering the cavity at the point $\xi$ with the velocity $-\vec v_0$. Note that $\Delta t = 2\pi/(\om m)$ is the minimal time period between two identical positions of the rotating set $B_m$. Then the momentum imparted to $B_m$ by the particles of the flow during the time interval $\Delta t$ equals
\begin{equation}\label{momentum}
2r\rho v \Delta t \int_0^1\!\!\int_{-\pi/2}^{\pi/2} (-\vec v_0 - \vec v^+_{(m)}(\xi, \vphi))\, \frac 12\, \cos\vphi\, d\vphi\, d\xi,
\end{equation}

Consider the frame of reference $\tilde O \tilde x_1 \tilde x_2$ having the center at the midpoint of the cavity opening $I$, the axis $\tilde O \tilde x_1$ being parallel to $I$, and $\tilde O \tilde x_2$, codirectional with $\vec n_\xi$. That is, the frame of reference rotates jointly with the segment $I$. The change of variables from $\vec x = (x_1;\, x_2)^T$ to $\vec{\tilde x} = (\tilde x_1;\, \tilde x_2)^T$ and  the reverse one are given by
$$
\vec{\tilde x} = A_{-\om t} \vec x - r\cos(\pi/m)\, e_{\pi/2} \quad \text{and} \quad \vec{x} = A_{\om t} \vec{\tilde x} + r\cos(\pi/m)\, e_{\pi/2+\om t},
$$
where $A_\phi = \left( \!\! \begin{array}{cr} \cos\phi &
-\sin\phi\\ \sin\phi & \cos\phi \end{array} \!\!
\right)$ and $\vec{e}_\phi =  \left[ \!\! \begin{array}{c} \cos\phi\\
\sin\phi \end{array} \!\! \right]$.

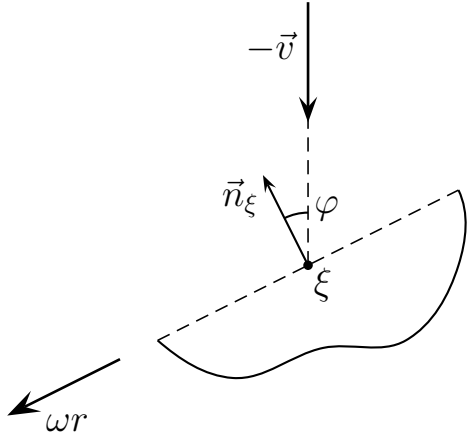
\begin{figure}[h]
\begin{picture}(0,170)
 \rput(5.6,1.4){
 \psline[linewidth=0.6pt,linestyle=dashed](0,0)(4,2)
 \psdot(2,1)
 \rput(2.2,0.75){\large $\xi$}
 \psecurve(-0.7,1)(0,0)(1,-0.5)(2.2,-0.1)(3.3,0)(4,1)(4,2)(3.5,2.4)
      \psline[linewidth=1pt,arrows=->,arrowscale=2.5](-0.5,-0.25)(-2,-1)
      \rput(-1.2,-1.1){\large $\om r$}
    \psline[linewidth=0.8pt,arrows=->,arrowscale=1.5](2,1)(1.4,2.2)
      \rput(1.1,1.9){\large $\vec n_\xi$}
  \psline[linewidth=0.6pt,linestyle=dashed](2,1)(2,3)
   \psarc[linewidth=0.8pt](2,1){0.7}{90}{116.57}
   \rput(2.25,1.75){\large $\vphi$}
   \psline[linewidth=1pt,arrows=->,arrowscale=2](2,4.5)(2,2.9)
   \rput(1.5,3.9){\large $-\vec v$}
   }
\end{picture}
\label{fig particle falling on a cavity} \caption{A particle falling on a cavity.}
\end{figure}
 \vspace{7mm}


Suppose now that $\vec x(t)$ and $\vec{\tilde{x}}(t)$ are  the  coordinates of a moving point in the initial and rotating frames of reference, respectively, and let $\vec v = \left( v_1;\, v_2 \right)^T = d\vec x/dt$ and $\vec{\tilde{v}} = \left( \tilde{v}_1;\, \tilde{v}_2 \right)^T = d\vec{\tilde{x}}/dt$. Then
\begin{equation}\label{transformation formulas}
\vec{\tilde v} = A_{-\om t} \vec{v} - \om A_{\pi/2-\om t} \vec x \ \ \ \text{and}
\ \ \ \vec{v} = A_{\om t} \vec{\tilde v} + \om A_{\pi/2+\om t} \vec{\tilde x}- \om r\cos({\pi}/{m})\, \vec{e}_{\om t}.
\end{equation}

We apply formulas (\ref{transformation formulas}) to the velocity of the particle at the two moments of its intersection with $I$. At the first moment, it holds $\om t = \vphi$ and $\vec x = re_{\pi/2+\vphi} + o(1)$ as $m \to \infty$. (Here and in what follows, the estimates $o(1)$ are not necessarily uniform with respect to $\xi$ and $\vphi$.) Then the incidence velocity $-\vec{v}_0$ takes the form
 \beq\label{v tilde}
-\vec{\tilde v}_0 = v \left( \lam - \sin\vphi;\ -\cos\vphi \right)^T + o(1) = -v\rrr \left( -\sin x;\,\cos x \right)^T + o(1),
 \eeq
where $\lam = \om r/v$,
 \beq\label{r and x}
\rrr = \rrr(\vphi, \lam) = \sqrt{\lam^2 - 2\lam\sin\vphi + 1} \ \
\text{and} \ \ x = x(\vphi, \lam) = \arcsin \frac{\lam -
\sin\vphi}{\rrr(\vphi, \lam)}.
 \eeq
As $m \to \infty$, the time spent by the particle in the cavity tends to zero, therefore the rotating frame of reference  can be considered ``approximately inertial'' during that time, the velocity at the second point of intersection given by
$$
\vec{\tilde v}^+ = v\rrr \left( -\sin y;\, \cos y \right)^T + o(1),\, \ \
\text{where}\, \ \ y = y(\xi, \vphi, \lam) = \vphi^+(\xi, x(\vphi, \lam)).
$$
(Here $\xi^+(\xi, \vphi)$,\, $\vphi^+(\xi, \vphi)$ denote the mapping generated by the cavity; see subsection 2.1.) Applying the second  formula in (\ref{transformation formulas}) and taking into account that $\vec{\tilde x} = o(1)$ and $\om t = \vphi + o(1)$, we find the velocity in the initial frame of reference
$$
\vec{v}^+ = \vec{v}_{(m)}^+(\xi, \vphi, \lam) = v\rrr A_\vphi (-\sin y;\, \cos y)^T - v\lam\, e_\vphi + o(1) = v^+(\xi, \vphi, \lam) + o(1),
$$
where
\begin{equation}\label{v+}
v^+(\xi, \vphi, \lam) = v\left[ \!\! \begin{array}{r} -\rrr \sin(\vphi + y) - \lam \cos\vphi\\
\rrr \cos(\vphi + y) - \lam \sin\vphi \end{array} \!\! \right].
\end{equation} \\

Letting $m \to \infty$ in the  formula (\ref{momentum}) for the  imparted momentum and dividing it by $\Delta t$, we get the following formula for the force of resistance acting on the disk
 \beq\label{integral R}
\vec{R} = \left[ \!\! \begin{array}{c} R_T\\ R_L \end{array} \!\!\! \right] =  r\rho v
\int_0^1\!\!\int_{-\pi/2}^{\pi/2} (-\vec{v}_0 - \vec{v}^+(\xi, \vphi, \lam))\, \cos\vphi\, d\xi\,
d\vphi.
 \eeq

The angular momentum transmitted to  $B_m$ by  an individual particle equals $ r v \rrr(\sin x + \sin y) + o(1)$  times the mass of the particle. Summing the angular momenta up over all incident particles and passing to the limit $m \to \infty$, one finds the moment of the resistance force acting on the disk,
  \beq\label{integral RI}
R_I\, =\, r^2\rho v
\int_0^1\!\!\int_{-\pi/2}^{\pi/2}
v\rrr(\vphi, \lam)\, (\sin x(\vphi, \lam) + \sin y(\xi, \vphi,
\lam))\, \cos\vphi\, d\xi\, d\vphi.
 \eeq

Below we will see that resistance and moment of resistance can be written down in the form
\begin{equation}\label{resistance}
\vec R = \frac 83\, r\rho v^2 \cdot \vec R[\nu_B, \lam],
\end{equation}
\begin{equation}\label{moment of resistance}
R_I = \frac 83\, r^2\rho v^2 \cdot R_I[\nu_B, \lam].
\end{equation}
The dimensionless values $\vec R[\nu, \lam] = (R_T[\nu, \lam];\, R_L[\nu, \lam])^T$ and $R_I[\nu, \lam]$ will be determined for several values of $\nu$ and $\lam$.

Let us consider separately the following three cases: $\lam = 1$,\, $0 < \lam < 1$, and $\lam > 1$, where the details of calculation are different.

\subsubsection{$\mathbf{\lam = 1}$}

In this case  we have $x = x(\vphi, 1) = \arcsin\sqrt{(1 - \sin\vphi)/2} = \pi/4
-\vphi/2$, and  so, the function $\vphi \mapsto x(\vphi, 1)$ is a bijection between the intervals $\interval$ and $[0,\, \pi/2]$. Further, one has $\rrr = \rrr(\vphi, 1) = \sqrt{2(1 - \sin\vphi)} = 2\sin x,\, \ \cos\vphi =
\sin 2x$, and we obtain from (\ref{v+})
$$
\vec v^+\, =\, v\left[ \!\! \begin{array}{r} -2\sin x\, \cos(2x - y) - \sin 2x\\
2\sin x\, \sin(2x - y) - \cos 2x \end{array} \!\! \right],
$$
wherefrom
$$
-\vec v_0 - \vec v^+\, =\, 2v\sin x \left[ \!\! \begin{array}{r} \cos(2x - y) + \cos x\\
-\sin(2x - y) - \sin x \end{array} \!\! \right].
$$
Making the change of variables $\{ \xi, \vphi \} \to \{ \xi, x \}$ in the integral  in (\ref{integral R}) and using (\ref{resistance}), one gets
$$
\vec R[\nu_B, 1]\, =\, 3\ \int_0^1\!\!\int_{0}^{\pi/2} \sin^2 x \left[ \!\!\! \begin{array}{r} \cos(2x - y) + \cos x\\
-\sin(2x - y) - \sin x \end{array} \!\! \right]\, \cos x\, d\xi\, dx.
$$
In this integral, $y$ is the function of $\xi$ and $x$,\, $y = \vphi^+(\xi, x)$. Changing the variables once again, $\{ \xi, x \} \to \{ x, y \}$, and taking into account that $\cos x\, d\xi\, dx = d\nu_B(x,y)$, we obtain
 \beq\label{R 1}
\vec R[\nu_B, 1] = 
3\int\!\!\!\int\limits_{\hspace*{-2mm}\square\hspace{-0.5mm}\blacksquare}
\sin^2 x \left[ \!\!
\begin{array}{r} \cos(2x-y) + \cos x\\ -\sin(2x-y) - \sin x
\end{array} \!\!
\right] d\nu_B(x,y)\,.
 \eeq
Here the symbol $\square\hspace{-0.5mm}\blacksquare$ stays for the rectangle $x \in [0,\, \pi/2], \ y \in \interval$.

The moment of the resistance force is calculated analogously, resulting in
$$
R_I[\nu_B, 1] = -3 \int_0^1\!\!\int_{0}^{\pi/2} \sin^2 x \left(\sin x + \sin y \right) \cos x\, d\xi\,
dx\, = \hspace*{40mm}
$$
 \beq\label{R_I 1}
\hspace*{40mm} =\, -3
\int\!\!\!\int\limits_{\hspace*{-2mm}\square\hspace{-0.5mm}\blacksquare}
\sin^2 x \left( \sin x + \sin y \right) d\nu_B(x,y).
 \eeq

\subsubsection{$\mathbf{0 < \lam < 1}$}

The second relation in (\ref{r and x}) implies that for a fixed value of $\lam$, $x = x(\vphi, \lam)$ is a monotone decreasing function of $\vphi$ that varies from $\pi/2$ to $-\pi/2$ as $\vphi$ changes from $-\pi/2$ to $\pi/2$. From formula (\ref{v
tilde}) and the first relation in (\ref{r and x}) we have
$$
\sin\vphi = \lambda \cos^2 x - \sin x \sqrt{1 - \lambda^2 \cos^2 x}
, \ \ \cos\vphi = \cos x\, (\lambda \sin x + \sqrt{1 - \lambda^2
\cos^2 x}),
$$
$$
\rrr = \lambda\sin x + \sqrt{1 - \lambda^2 \cos^2 x}.
$$
Denote
\begin{equation}\label{zeta}
\zeta = \zeta(x, \lam) = \arcsin \sqrt{1 - \lambda^2 \cos^2 x};
\end{equation}
then one has
$$
\cos\zeta = \lambda\cos x, \ \ \  x + \zeta = \pi/2 - \vphi, \ \ \
\zeta \in [\arccos\lambda,\, \pi/2],
$$
and taking into account (\ref{v+}), we get
$$
-\vec v_0 - \vec v^+\, =\, v(\lambda\sin x + \sin\zeta) \cdot
2\cos\frac{x-y}{2} \left[ \!\!
\begin{array}{r} \cos\left(\zeta + \frac{x-y}{2} \right)\\ -\sin\left(\zeta + \frac{x-y}{2}
\right)
\end{array}
\!\! \right],
$$
$$
\frac{\cos\vphi}{\cos x}\, =\, \lambda\sin x + \sin\zeta\, =\, \rrr,
$$
$$
\frac{d\vphi}{dx}\ =\ -1 - \frac{d\zeta}{dx}\ =\ -\,
\frac{\lambda\sin x + \sin\zeta}{\sin\zeta}\,.
$$
\\

Using the obtained formulas, making the change of variables $\{ \xi, \vphi \} \to \{ \xi, x \}$ in the integral (\ref{integral R}), and taking into account (\ref{resistance}), one gets
$$
\vec R[\nu_B, \lam] = \frac 34 \int_0^1\!\!\!\int_{-\pi/2}^{\pi/2}
\frac{(\lambda \sin x + \sin\zeta)^3}{\sin\zeta}\, \cos\frac{x-y}{2}\!
\left[ \hspace{-1.5mm}
\begin{array}{r} \cos\left(\zeta + \frac{x-y}{2} \right)\\ -\sin\left(\zeta + \frac{x-y}{2}
\right)
\end{array} \hspace{-1.5mm}
\right]\! \cos x\, d\xi\, dx.
$$
\\

Finally,  the change of variables $\{ \xi, x \} \to \{ x, y \}$  results in
\beq\label{R lam<1}
\vec R[\nu_B, \lam] =
\frac 34\! \int\!\!\!\!\int\limits_{\hspace*{-2mm}\Box}\!
\frac{(\lambda \sin x + \sin\zeta)^3}{\sin\zeta} \cos\frac{x-y}{2}\!
\left[ \hspace{-2mm}
\begin{array}{r} \cos\left(\zeta + \frac{x-y}{2} \right)\\ -\sin\left(\zeta + \frac{x-y}{2}
\right)
\end{array} \hspace{-2mm}
\right]\! d\nu_B(x,y).
 \eeq
Recall that the symbol  $\Box$ denotes the square $\interval \times \interval$
and $\zeta = \zeta(x, \lam)$.

In a similar way, from (\ref{integral RI}) one gets
$$
R_I\, =\, -\frac 38 \int_0^1\!\!\int_{-\pi/2}^{\pi/2} \frac{(\lambda \sin x +
\sin\zeta)^3}{\sin\zeta}\, (\sin x + \sin y)\, \cos x\, d\xi\, dx;
$$
wherefrom
 \beq\label{R_I lam<1}
R_I[\nu_B, \lam]\, =\, -\frac 38
\int\!\!\!\int\limits_{\hspace*{-2mm}\Box} \frac{(\lambda \sin x +
\sin\zeta)^3}{\sin\zeta} \left( \sin x + \sin y \right)\,
d\nu_B(x,y).
 \eeq

Formulas (\ref{R 1}) and (\ref{R_I 1}) are the particular cases of (\ref{R lam<1}) and (\ref{R_I lam<1}) for
$\lam = 1$. This can be easily verified taking into account that $\zeta(x, 1) = |x|$.

\subsubsection{$\mathbf{\lam > 1}$}

In this case, $x = x(\vphi, \lam)$ (\ref{r and x}) is not injection anymore. When $\vphi$ varies from $-\pi/2$
to $\vphi_0 = \vphi_0(\lam) := \arcsin\frac{1}{\lam}$, the value of $x$ monotonically decreases from $\pi/2$ to
$x_0 = x_0(\lam) := \arccos\frac{1}{\lam}$, and when $\vphi$ varies from $\vphi_0$ to $\pi/2$, $\ x$ monotonically increases from $x_0$ to $\pi/2$. Denote by $\vphi_- := \vphi_-(x, \lam)$ and $\vphi_+:=\vphi_+(x, \lam)$ the functions inverse to $x(\vphi, \lam)$ on the intervals $[-\pi/2,\, \vphi_0]$ and $[\vphi_0,\, \pi/2]$, respectively. Then one has
$$
\sin\vphi_\pm\, =\, \lambda \cos^2 x\, \pm\, \sin x \sqrt{1 -
\lambda^2 \cos^2 x}.
$$
Here and in what follows the signs ``$+$'' and ``$-$'' are related to the functions $\vphi_+$ and $\vphi_-$, respectively. The values $\vphi_+$,\, $\vphi_-$, and $\zeta = \zeta(x, \lam)$ (\ref{zeta}) satisfy the relations
$$
\pi/2 - \vphi_+\, =\, x - \zeta, \ \ \ \ \ \pi/2 - \vphi_-\, =\, x +
\zeta.
$$
The function $\zeta$ is defined for $x \in [x_0,\, \pi/2]$ and monotonically increases from 0 to $\pi/2$, when $x$ changes in the interval $[x_0,\, \pi/2]$.

After some algebra one gets
$$
\frac{\cos\vphi_\pm}{\cos x}\ =\ \frac{\sin(x \mp \zeta)}{\cos x}\,
=\, \lambda \sin x \mp \sin\zeta;
$$
$$
\pm \frac{d\vphi_\pm}{dx}\ =\ \frac{d\zeta}{dx}\, \mp\, 1\ =\
\frac{\lambda\sin x \mp \sin\zeta}{\sin\zeta};
$$
$$
\rrr_\pm\, =\, \lambda\sin x \mp \sin\zeta;
$$
$$
-\vec v_0 - \vec v^+_±\, =\, v\left( \lambda\sin x \mp \sin\zeta
\right) \cdot 2\cos\frac{x-y}{2}
\left[\! \hspace{-1mm}
\begin{array}{r} \cos\left(\frac{x-y}{2} \mp \zeta \right)\\
-\sin\left(\frac{x-y}{2} \mp \zeta \right)
\end{array} \hspace{-1mm}
\! \right].
$$
\vspace{3mm}

Here the shorthand notation $\rrr_\pm = \rrr(\vphi_\pm(x, \lam), \lam)$,\, $\vec
v^+_± = \vec v^+(\xi, \vphi_\pm(x, \lam), \lam)$,\, $y = y(\xi,
\vphi_\pm(x, \lam), \lam) = \vphi_i^+(\xi, x)$ is used.

The resistance force  takes the form $\vec R[\nu_B, \lam] = \vec R_- + \vec R_+$, where
$$
\vec R_\pm\, =\, \frac 38 \int_0^1\!\!\int_{x_0}^{\pi/2} (-\vec v_0 - \vec v^+_±)\ \frac{\cos\vphi_\pm}{\cos x} \left( \! \pm
\frac{d\vphi_\pm}{dx} \! \right) \cos x\, d\xi\, dx\, =
$$
$$
=\, \frac 34\, \int_0^1\!\!\int_{x_0}^{\pi/2} \frac{(\lambda \sin x
\mp \sin\zeta)^3}{\sin\zeta}\ \cos\frac{x-y}{2} \left[ \hspace{-1.5mm}
\begin{array}{c} \cos\left(\frac{x-y}{2} \mp \zeta \right)\\
-\sin\left(\frac{x-y}{2} \mp \zeta \right)
\end{array} \hspace{-1.9mm}
\right] \cos x\, d\xi\, dx.
$$ \\

Summing the integrals above and making the change of variables, one obtains
$$
\vec R[\nu_B, \lam] = \frac 32
\int\!\!\!\int\limits_{\hspace*{-2mm}\square\hspace{-0.5mm}\blacksquare}
\frac{\cos\frac{x-y}{2}}{\sin\zeta}\ \left\{ (\lam^3 \sin^3 x +
3\lam\sin x \sin^2 \zeta) \cos\zeta
\left[ \hspace{-1mm} \begin{array}{c} \cos\frac{x-y}{2}\\
-\sin\frac{x-y}{2}
\end{array} \hspace{-1.5mm}
\right] \right. -
$$
 \beq\label{R i lam>1}
- \left. (3\lam^2 \sin^2 x \sin\zeta + \sin^3 \zeta) \sin\zeta
\left[ \hspace{-1mm}
\begin{array}{c}
\sin\frac{x-y}{2}\\
\cos\frac{x-y}{2}
\end{array} \hspace{-1.5mm}
\right] \right\} d\nu_B(x,y).
 \eeq \\

Here the symbol $\square\hspace{-0.5mm}\blacksquare$ stands for the rectangle $[x_0,\, \pi/2] \times \interval$.

The moment of the resistance force is calculated analogously. One has $R_I[\nu_B, \lam] = R_{I-} + R_{I+}$, where
$$
R_{I\pm}\, =\, -\frac 38 \int_0^1\!\!\int_{x_0}^{\pi/2} \rrr_\pm\
\frac{\cos\vphi_\pm}{\cos x} \left( \! \pm \frac{d\vphi_\pm}{dx} \!
\right) (\sin x + \sin y) \cos x\, d\xi\, dx\, =
$$
$$
=\, -\frac 38 \int_0^1\!\!\int_{x_0}^{\pi/2} \frac{(\lambda \sin x
\mp \sin\zeta)^3}{\sin\zeta} \left( \sin x + \sin y \right)
 \cos x\, d\xi\, dx.
$$
Therefore
$$
R_{I}[\nu_B, \lam] = -\frac 34 \int_0^1\!\!\int_{x_0}^{\pi/2} \frac{\lambda^3 \sin^3 x + 3\lam \sin x \sin^2 \zeta}{\sin\zeta} \left( \sin x + \sin y \right) \cos x\, d\xi\, dx.
$$
Making the change of variables, we have
 \beq\label{R I lam>1}
R_I[\nu_B, \lam] = -\frac 34
\int\!\!\!\int\limits_{\hspace*{-2mm}\square\hspace{-0.5mm}\blacksquare}
\frac{\lambda^3 \sin^3 x + 3\lam \sin x \sin^2 \zeta}{\sin\zeta}
\left( \sin x + \sin y \right) d\nu_B(x,y).
 \eeq

\section{Magnus effect}

 We are concerned here with calculation of the ``normalized'' resistance force $\vec R[\nu, \lam]$ and its moment $R_I[\nu, \lam]$ for some special cases of $\nu$ and $\lam$, as well as determination of the two-dimensional set $\{ \vec R[\nu, \lam]: \ \nu \in \Pigam^{\text{symm}} \}$ for several values of $\lam$.

In the previous section, the formulas for $\vec R[\nu, \lam]$ and $R_I[\nu, \lam]$ were derived,
 \beq\label{functionals1}
\vec R[\nu, \lam] = \left[\!\! \begin{array}{c} R_T[\nu, \lam]\\ R_L[\nu, \lam] \end{array} \!\!\right] = \int\!\!\!\int\limits_{\hspace*{-2mm}\Box}
\vec c(x, y, \lam)\, d\nu(x,y)
 \eeq
and
 \begin{equation*}\label{functionals2}
 R_I[\nu, \lam] = \int\!\!\!\int\limits_{\hspace*{-2mm}\Box}
c_I(x, y, \lam)\, d\nu(x,y),
 \end{equation*}
with $\vec c$ and $c_I$ given by  the relations  below (\ref{c lam<1})--(\ref{c_I lam>1}).
Recall that $\lam = \om r/v$ is the relative angular velocity and the symbol $\Box$ denotes the square $\interval \times \interval$. In the formulas below, the value $\zeta = \zeta(x, \lam)$ is given by (\ref{zeta}), $x_0 = x_0(\lam) = \arccos\frac{1}{\lam}$, and $\chi$ stands for the characteristic function.\\

(a) If $0 < \lam \le 1$ then
 \beq\label{c lam<1}
\vec c(x, y, \lam)\ =\ \frac 34\ \frac{(\lambda \sin x +
\sin\zeta)^3}{\sin\zeta}\, \cos\frac{x-y}{2} \left[ \hspace{-1.5mm}
\begin{array}{r} \cos\left(\zeta + \frac{x-y}{2} \right)\\
-\sin\left(\zeta + \frac{x-y}{2} \right)
\end{array} \hspace{-1.5mm}
\right],
 \eeq
 \beq\label{c_I lam<1}
c_I(x, y, \lam)\ =\ -\frac 38\ \frac{(\lambda \sin x +
\sin\zeta)^3}{\sin\zeta} \left( \sin x + \sin y \right),
 \eeq
and in particular,
 \beq\label{c lam=1}
\vec c(x, y, 1) = 3\sin^2 x\! \left[ \!\!\!
\begin{array}{r} \cos(2x-y) + \cos x\\
-\sin(2x-y) - \sin x \end{array} \!\!\!
\right] \chi_{x \ge 0}(x,y),
 \eeq
 \beq\label{c_I lam=1}
c_I(x, y, 1) = -3\sin^2 x \left( \sin x + \sin y \right) \chi_{x \ge
0}(x,y).
 \eeq
In the limiting case $\lam \to 0^+$ one has
 \beq\label{c lam=0}
\vec c(x,y,\lam)\, =\, -\frac{3}{8} \left[ \!\! \begin{array}{c} \sin(x-y)\\ 1 + \cos(x-y) \end{array} \!\! \right] + O(\lam),
 \eeq
 \vspace{-3mm}
$$
c_I(x,y)\, =\,  \frac{9\lambda}{8}\, \sin x \left( \sin x + \sin y \right) + O(\lam^2). \quad \quad \quad \quad \quad
$$

(b) In the case $\lam > 1$ we have
$$
\vec c(x, y, \lam)\ =\ \frac 32\ \frac{\cos\frac{x-y}{2}}{\sin\zeta}
\left\{ (\lam^3 \sin^3 x + 3\lam\sin x \sin^2 \zeta) \cos\zeta
\left[ \hspace{-1mm} \begin{array}{c} \cos\frac{x-y}{2}\\
-\sin\frac{x-y}{2}
\end{array} \hspace{-1mm}
\right] \right. -
$$
 \beq\label{c lam>1}
\quad  \quad - \left. (3\lam^2 \sin^2 x \sin\zeta + \sin^3 \zeta) \sin\zeta
\left[ \hspace{-1mm}
\begin{array}{c}
\sin\frac{x-y}{2}\\
\cos\frac{x-y}{2}
\end{array} \hspace{-1.5mm}
\right] \right\} \chi_{x \ge x_0}(x,y),
 \eeq \\
  \beq\label{c_I lam>1}
c_I(x, y, \lam)\, =\, -\frac 34\ \frac{\lambda^3 \sin^3 x + 3\lam \sin x
\sin^2 \zeta}{\sin\zeta} \left( \sin x + \sin y \right) \chi_{x \ge
x_0}(x,y).
 \eeq

 \subsection{Special cases of rough circles}

Consider several examples of rough circles. We present here only the final expressions for the normalized forces and their moments  calculated by  formulas (\ref{c lam<1})--(\ref{c_I lam>1}); the calculation details are omitted.

\begin{enumerate}

\item {\it Circle}.~ Consider the measure $\nu_0 \in \Pigam^{\text{symm}}$ corresponding to the  ordinary (without roughness) circle. The measure $\nu_0$ is supported on the diagonal $y = -x$; the density of the measure is $d\nu_0(x, y) = \cos x \cdot \del(x + y)$. One has $R_T[\nu_0, \lam] = R_I[\nu_0, \lam] = 0$ and $R_L[\nu_0, \lam] = -1$. Thus, as one could expect, the resistance does not depend on the angular velocity and is collinear to the body's velocity. There is no Magnus effect in this case.

\item {\it Retroreflector}.~ There exists a unique measure $\nu_\star \in \Pigam^{\text{symm}}$ supported on the diagonal $x = y$; its density equals $d\nu_\star(x, y) = \cos x \cdot \del(x - y)$. We believe that there is no rough circle generating this measure; but, according to Theorem \ref{teor2}, there exist rough circles {\em approximating} it; that is, there exists a sequence of rough circles $B_n$ such that $\nu_{B_n}$ weakly converge to $\nu_\star$ (see also \cite{Nonlinearity}). One has $R_T[\nu_\star, \lam] = {3\pi}\lam/{8}$,\, $R_L[\nu_\star, \lam] = -3/2$ and $R_I[\nu_\star, \lam] = -3\lambda/2$. Thus, the longitudinal component of the resistance force does not depend on the angular velocity $\lam$, while the transversal component and the moment of this force are proportional to $\lam$.

\item {\it Rectangular cavity}.~ The cavity is rectangular, with the width being much smaller that the depth, ({\it width})/({\it height of the rectangle})$\ = \ve$ for $\ve>0$ small (see Fig. 5\,(a)). We also suppose that the common length of the convex part of the boundary $\pl B$ is very small; say, equal to $\ve$. Then $\nu_{B} = \nu_{\text{rect}} + o(1)$, where $\nu_{\text{rect}} = (\nu_0 + \nu_\star)/2$ and $ o(1)$ denotes a measure that weakly converges to zero as $\ve \to 0^+$. One can easily calculate that $R_T[\nu_{\text{rect}}, \lam] = {3\pi}\lam/{16}$,\, $R_L[\nu_{\text{rect}}, \lam] = -1.25$ and $R_I[\nu_{\text{rect}}, \lam] = -3\lambda/4$.

\begin{figure}[h]
\vspace{3mm}
\begin{picture}(0,150)
\rput(3.7,2.5){
\scalebox{0.4}{
\rput{-150}
{\psline[linestyle=dashed,linewidth=0.4pt](1.6077,6)(0,0)(-1.6077,6)
\psline[linestyle=dotted,linewidth=1.6pt](0.55,6)(-0.55,6)
\psline[linewidth=2.4pt](1.6077,6)(0.55,6)(0.55,3.6)(-0.55,3.6)(-0.55,6)(-1.6077,6)}
\rput{-120}
{\psline[linestyle=dashed,linewidth=0.4pt](1.6077,6)(0,0)(-1.6077,6)
\psline[linestyle=dotted,linewidth=1.6pt](0.55,6)(-0.55,6)
\psline[linewidth=2.4pt](1.6077,6)(0.55,6)(0.55,3.6)(-0.55,3.6)(-0.55,6)(-1.6077,6)}
\rput{-90}
{\psline[linestyle=dashed,linewidth=0.4pt](1.6077,6)(0,0)(-1.6077,6)
\psline[linestyle=dotted,linewidth=1.6pt](0.55,6)(-0.55,6)
\psline[linewidth=2.4pt](1.6077,6)(0.55,6)(0.55,3.6)(-0.55,3.6)(-0.55,6)(-1.6077,6)}
\rput{-60}
{\psline[linestyle=dashed,linewidth=0.4pt](1.6077,6)(0,0)(-1.6077,6)
\psline[linestyle=dotted,linewidth=1.6pt](0.55,6)(-0.55,6)
\psline[linewidth=2.4pt](1.6077,6)(0.55,6)(0.55,3.6)(-0.55,3.6)(-0.55,6)(-1.6077,6)}
\rput{-30}
{\psline[linestyle=dashed,linewidth=0.4pt](1.6077,6)(0,0)(-1.6077,6)
\psline[linestyle=dotted,linewidth=1.6pt](0.55,6)(-0.55,6)
\psline[linewidth=2.4pt](1.6077,6)(0.55,6)(0.55,3.6)(-0.55,3.6)(-0.55,6)(-1.6077,6)}
\rput{180}
{\psline[linestyle=dashed,linewidth=0.4pt](1.6077,6)(0,0)(-1.6077,6)
\psline[linestyle=dotted,linewidth=1.6pt](0.55,6)(-0.55,6)
\psline[linewidth=2.4pt](1.6077,6)(0.55,6)(0.55,3.6)(-0.55,3.6)(-0.55,6)(-1.6077,6)}
\rput{150}
{\psline[linestyle=dashed,linewidth=0.4pt](1.6077,6)(0,0)(-1.6077,6)
\psline[linestyle=dotted,linewidth=1.6pt](0.55,6)(-0.55,6)
\psline[linewidth=2.4pt](1.6077,6)(0.55,6)(0.55,3.6)(-0.55,3.6)(-0.55,6)(-1.6077,6)}
\rput{120}
{\psline[linestyle=dashed,linewidth=0.4pt](1.6077,6)(0,0)(-1.6077,6)
\psline[linestyle=dotted,linewidth=1.6pt](0.55,6)(-0.55,6)
\psline[linewidth=2.4pt](1.6077,6)(0.55,6)(0.55,3.6)(-0.55,3.6)(-0.55,6)(-1.6077,6)}
\rput{90}
{\psline[linestyle=dashed,linewidth=0.4pt](1.6077,6)(0,0)(-1.6077,6)
\psline[linestyle=dotted,linewidth=1.6pt](0.55,6)(-0.55,6)
\psline[linewidth=2.4pt](1.6077,6)(0.55,6)(0.55,3.6)(-0.55,3.6)(-0.55,6)(-1.6077,6)}
\rput{60}
{\psline[linestyle=dashed,linewidth=0.4pt](1.6077,6)(0,0)(-1.6077,6)
\psline[linestyle=dotted,linewidth=1.6pt](0.55,6)(-0.55,6)
\psline[linewidth=2.4pt](1.6077,6)(0.55,6)(0.55,3.6)(-0.55,3.6)(-0.55,6)(-1.6077,6)}
\rput{30}
{\psline[linestyle=dashed,linewidth=0.4pt](1.6077,6)(0,0)(-1.6077,6)
\psline[linestyle=dotted,linewidth=1.6pt](0.55,6)(-0.55,6)
\psline[linewidth=2.4pt](1.6077,6)(0.55,6)(0.55,3.6)(-0.55,3.6)(-0.55,6)(-1.6077,6)}
\psline[linestyle=dashed,linewidth=0.4pt](1.6077,6)(0,0)(-1.6077,6)
\psline[linestyle=dotted,linewidth=1.6pt](0.55,6)(-0.55,6)
\psline[linewidth=2.4pt](1.6077,6)(0.55,6)(0.55,3.6)(-0.55,3.6)(-0.55,6)(-1.6077,6)
}
 \vspace{32mm}
   }
\rput(11,2.5){
\scalebox{0.8}{
\rput{-150}
{\psline[linestyle=dashed,linewidth=0.4pt](0.80385,3)(0,0)(-0.80385,3)
\psline[linewidth=1.2pt](0.80385,3)(0,2.196)(-0.80385,3)}
\rput{-120}
{\psline[linestyle=dashed,linewidth=0.4pt](0.80385,3)(0,0)(-0.80385,3)
\psline[linewidth=1.2pt](0.80385,3)(0,2.196)(-0.80385,3)}
\rput{-90}
{\psline[linestyle=dashed,linewidth=0.4pt](0.80385,3)(0,0)(-0.80385,3)
\psline[linewidth=1.2pt](0.80385,3)(0,2.196)(-0.80385,3)}
\rput{-60}
{\psline[linestyle=dashed,linewidth=0.4pt](0.80385,3)(0,0)(-0.80385,3)
\psline[linewidth=1.2pt](0.80385,3)(0,2.196)(-0.80385,3)}
\rput{-30}
{\psline[linestyle=dashed,linewidth=0.4pt](0.80385,3)(0,0)(-0.80385,3)
\psline[linewidth=1.2pt](0.80385,3)(0,2.196)(-0.80385,3)}
\rput{180}
{\psline[linestyle=dashed,linewidth=0.4pt](0.80385,3)(0,0)(-0.80385,3)
\psline[linewidth=1.2pt](0.80385,3)(0,2.196)(-0.80385,3)}
\rput{150}
{\psline[linestyle=dashed,linewidth=0.4pt](0.80385,3)(0,0)(-0.80385,3)
\psline[linewidth=1.2pt](0.80385,3)(0,2.196)(-0.80385,3)}
\rput{120}
{\psline[linestyle=dashed,linewidth=0.4pt](0.80385,3)(0,0)(-0.80385,3)
\psline[linewidth=1.2pt](0.80385,3)(0,2.196)(-0.80385,3)}
\rput{90}
{\psline[linestyle=dashed,linewidth=0.4pt](0.80385,3)(0,0)(-0.80385,3)
\psline[linewidth=1.2pt](0.80385,3)(0,2.196)(-0.80385,3)}
\rput{60}
{\psline[linestyle=dashed,linewidth=0.4pt](0.80385,3)(0,0)(-0.80385,3)
\psline[linewidth=1.2pt](0.80385,3)(0,2.196)(-0.80385,3)}
\rput{30}
{\psline[linestyle=dashed,linewidth=0.4pt](0.80385,3)(0,0)(-0.80385,3)
\psline[linewidth=1.2pt](0.80385,3)(0,2.196)(-0.80385,3)}
\psline[linestyle=dashed,linewidth=0.4pt](0.80385,3)(0,0)(-0.80385,3)
\psline[linewidth=1.2pt](0.80385,3)(0,2.196)(-0.80385,3)
}
}
\end{picture}
\label{fig rectangular triangular cavity}
\caption{(a)  A rough disk with  \hspace*{30mm} (b) A rough disk with\newline \hspace*{20mm} rectangular cavities. \hspace*{35mm} triangular cavities.}
\vspace{5mm}
\end{figure}

\item {\it Triangular cavity}.~ The cavity is a right isosceles triangle (see Fig. 5\,(b)). The measure $\nu_{B} =: \nu_{\triangledown}$ has the following support (looking like an inclined letter  {\bf H}):~ $\{ x + y = -\pi/2: \  x \in [-\pi/2,\, 0] \} \cup \{ y = x: \ x \in [-\pi/4,\, \pi/4] \} \cup \{ x + y = \pi/2: \ x \in [0,\, \pi/2] \}$. The density of this measure equals
$d\nu_{\triangledown}(x,y) = \cos x \cdot \left( \chi_{[-\frac{\pi}{2},-\frac{\pi}{4}]}(x)\, \del(x + y +
\frac \pi 2) + \chi_{[-\frac{\pi}{4},\frac{\pi}{4}]}(x)\cdot \right.$ $\cdot \del(x - y) + $
$\left. \chi_{[\frac{\pi}{4},\frac{\pi}{2}]}(x)\, \del(x + y - \frac
\pi 2) \right)\, +\, |\sin x| \cdot \left(
\chi_{[-\frac{\pi}{4},0]}(x)\, \del(x + y + \frac
\pi 2) - \right.$\\
$- \left. \chi_{[-\frac{\pi}{4},\frac{\pi}{4}]}(x)\, \del(x - y) +
\chi_{[0,\frac{\pi}{4}]}(x)\, \del(x + y - \frac \pi 2) \right)$.
One has $\vec R[\nu_\triangledown, 0^+] = (0; \ -\sqrt{2})^T$ and $R_I[\nu_\triangledown, 0^+] = 0$;~ $\vec R[\nu_\triangledown, 1] = (1/4 + 3\pi/16; \ 3\pi/16 - 2)^T$. The rest of the values are still unknown.

\item {\it Cavity realizing the product-measure}.~ Consider the measure $\nu_\otimes$ with the density
$d\nu_{\otimes}(x,y) = \frac 12\, \cos x \cos y\, dx\, dy$. Evidently in this case $\nu_\otimes \in \Pigam^{\text{symm}}$. The angles of incidence and of reflection are statistically independent; so to speak, at the moment  when the particle leaves the cavity, it completely ``forgets'' its initial velocity. Here we have $R_T[\nu_\otimes, \lam] = (10\lam + \lam^3)\, \pi/80$ for $0 < \lambda \le 1$; $\,R_L[\nu_{\otimes}, 1] = -3/4 - \pi/5 \approx -1.378$; and $R_I[\nu_{\otimes}, \lam] = -3\lam/4$ for any $\lam$. The remaining values are unknown.
\end{enumerate}

\subsection{Vector-valued Monge-Kantorovich problem}

Consider the following problem: determine the set of all possible resistance forces  that can act on a rough circle. The angular velocity of the circle is fixed, and the roughness varies. The force is scaled so that the resistance of the  ``ordinary circle'' equals $(0;\, -1)^T$. Thus, the problem can be reformulated in the following way: given $\lam$, find the two-dimensional set
 \beq\label{eM 2}
\mathcal{R}_\lam := \{ \vec R[\nu, \lam]: \ \nu \in \Pigam^{\text{symm}} \}.
 \eeq
It can be viewed as a restriction of the following more general problem: find the three-dimensional set
 \begin{equation*}\label{eM 1}
\{ (\vec R[\nu, \lam];\, R_I[\nu, \lam]) : \ \nu \in \Pigam^{\text{symm}} \}.
 \end{equation*}
The latter problem is more important, but also more time-consuming, and is {mainly} postponed to the future. The only exception is the case $\lam = 1$; on Fig.\,8 several ``level sets'' $\mathcal{R}_{1,c} = \{ \vec R[\nu, 1]: \ \nu \in \Pigam^{\text{symm}}, R_I[\nu, 1] = c \}$ are depicted, suggesting an idea how the corresponding three-dimensional set looks like. In this case $R_I[\nu, 1]$ varies between $-1.5$ and 0, and the level sets are found for 21 values $c = -1.5,\,  -1.425,\, -1.35, \ldots, -0.15,\,  -0.075,\, 0$.

Note that the functional $\vec R$, defined on the set  $\Pigam^\text{symm}$ by  formula (\ref{functionals1}), will not change if the integrand $\vec c$ is replaced with the symmetrized function $\vec c^{\, \text{symm}}(x,y,\lam) = \frac 12 (\vec c(x,y,\lam) + \vec c(y,x,\lam))$:
 \begin{equation*}\label{eM 0}
\vec R[\nu,\lam] = \int\!\!\!\int\limits_{\hspace*{-2mm}\Box} \vec c^{\, \text{symm}}(x,y,\lam)\, d\nu(x,y).
  \end{equation*}
Denote by $\Pigam$ the set of measures $\nu$ on the square $\Box$ that satisfy
the only condition (\ref{uslovie 1}), that is, the set of measures with both marginals equal to $\gam$.
For any $\nu \in \Pigam$  it holds $\int\!\!\!\int\limits_{\hspace*{-2mm}\Box}
  \vec c^{\, \text{symm}}\, d\nu = \int\!\!\!\int\limits_{\hspace*{-2mm}\Box}
  \vec c^{\, \text{symm}}\, d\nu^{\text{symm}}$, where $\nu^{\text{symm}} = \frac 12 (\nu + \pi_d^\# \nu) \in \Pigam^{\text{symm}}$.
It follows that
 \beq\label{eM 3}
\mathcal{R}_\lam = \Big\{ \int\!\!\!\!\int\limits_{\hspace*{-2mm}\Box}
 \vec c^{\,\text{symm}}(x,y,\lam)\, d\nu(x,y): \ \ \nu \in \Pigam \Big\}.
 \eeq

The problem of finding $\mathcal{R}_\lam$ in (\ref{eM 3}) is a vector-valued analogue of the Monge-Kantorovich problem. The difference consists in the fact that the cost function, and therefore the functional, are vector valued. The set $\mathcal{R}_\lam$ is convex, since it is the image of the convex set $\Pigam$ under a linear mapping.

Note that, due to formula (\ref{c lam=0}), $\vec c^{\, \text{symm}}(x,y,0^+) = \frac 38\, (1 + \cos(x-y)) (0;\, 1)^T$, therefore the problem of finding $\mathcal{R}_{0^+}$ amounts to minimizing and maximizing the integral $\frac 38 \displaystyle\int\!\!\!\int\limits_{\hspace*{-2mm}\Box} (1 + \cos(x-y))\, d\nu(x,y)$ over all $\nu \in \Pigam$. This special Monge-Kantorovich problem was solved in \cite{MatSb 04: mean resist}; the minimal and maximal values of the integral were found to be $0.9878...$ and $1.5$.

On Figures 6 and 7 we present numerical solutions of this problem for the values $\lam = 0.1,\ 0.3$ and 1, as well as the analytical solution for $\lam = 0^+$. The case of  larger $\lam$ requires more involved calculation and therefore is postponed to the future. The method of solution is the following: for $n$ equidistant vectors $\vec e_i$,\, $i = 1,\ldots,n$ on $S^1$, we find the solution of the Monge-Kantorovich problem $\inf \langle \vec R[\nu,\lam],\, \vec e_i \rangle =: r_i$. Here $\langle \cdot\,, \cdot \rangle$ stands for the scalar product. This problem is reduced to the  transport problem of Linear Programming and is solved numerically.\footnote{ All the computational tests were performed on a PC
Pentium IV, 2.0Ghz and 512 Mb RAM and using the optimization package Xpress-IVE, Version 1.19.00 with the modeler MOSEL.}

Next, the intersection of the half-planes $\langle \vec r,\, \vec e_i \rangle \ge r_i$ is built. It is a convex polygon approximating the required set $\mathcal{R}_\lam$, and the approximation accuracy increases as $n$ increases. The value $n = 100$ was used in our calculations.

 On Fig.\,6, the sets $\mathcal{R}_\lam$ are shown for $\lam = 0^+,\, 0.1,\, 0.3$, and 1. The set $\mathcal{R}_{0^+}$ is the vertical segment $\{ 0 \} \times [-1.5,\, -0.9878]$,\, $\mathcal{R}_{0.1}$ is the thin set with white interior, and $\mathcal{R}_{0.3}$ is the set with gray interior.  The largest set is $\mathcal{R}_{1}$.

On Fig.\,7, the same sets are shown in more detail. On Figures 7\,(b)--(d), additionally, we present the regions corresponding to all possible kinds of roughness related to triangular cavities (and to cavities formed by combinations of different triangles), with the angles being multiples of $5^0$. These regions are colored red. For $\lam = 0^+$, the corresponding region is the vertical interval $\{ 0 \} \times [-1.42,\, -1]$ marked by a (slightly shifted) dashed line on Fig.\,7(a). The points corresponding to special values of $\nu$ are indicated by  special signs: $\nu_0$ is marked by a circle, $\nu_\star$ is marked by a dot, $\nu_{\text{rect}}$, by a square, $\nu_\triangledown$, by a triangle, and $\nu_\otimes$, by a circumscribed cross.

The part of the set  $\mathcal{R}_\lam$ situated to the left of the vertical axis  is related to  resistance forces  producing the (direct) Magnus effect. The part of $\mathcal{R}_\lam$ situated to the right of this axis, is related to  resistance forces  that cause the reverse Magnus effect. We can see that the most part of the set (in the case $\lam = 1$, approximately 93.6$\%$ of the area) is situated to the right of the axis. This suggests that the {\it reverse} effect is more common phenomenon than the direct one.  Actually, although Theorem 2 guarantees existence of an everywhere dense subset of $\mathcal{R}_\lam$ generated by shapes of roughness, we never met a simple shape producing the {\it direct} Magnus effect (and thus corresponding to a point on the left of the vertical axis).


\begin{figure}[!h]
\begin{center}
\includegraphics[width=0.9\textwidth]{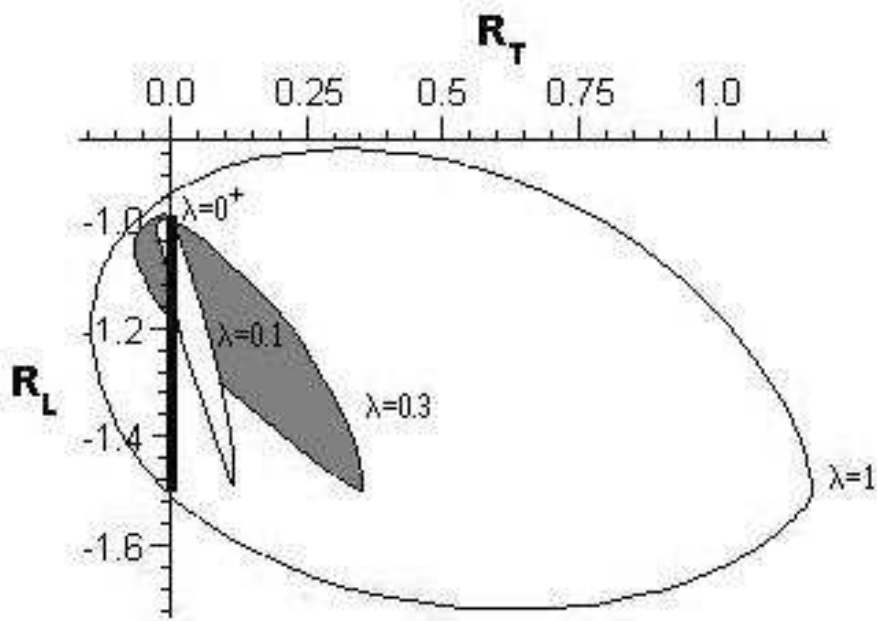}
\end{center}
\caption{The convex sets $\mathcal{R}_\lam$ with $\lam = 0^+,\, 0.1,\, 0.3$ and 1 are shown. The set $\mathcal{R}_{0^+}$ is the vertical segment with the endpoints $(0, - 0.9878...)$ and $(0, -1.5)$.}
\end{figure}

 \newpage

\begin{figure}[htb!]
\begin{center}
\subfigure[][$\lam = 0$]{\includegraphics[width=.40\textwidth]{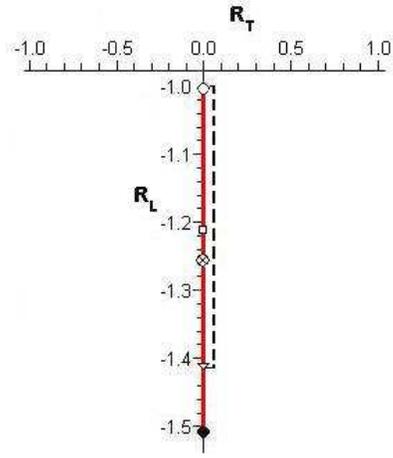}}
\hfill
\subfigure[][$\lam = 0.1$]{\includegraphics[width=.40\textwidth]{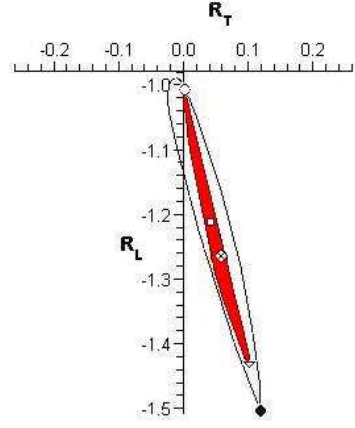}}
\\
\subfigure[$\lam = 0.3$]{\includegraphics[width=.40\textwidth]{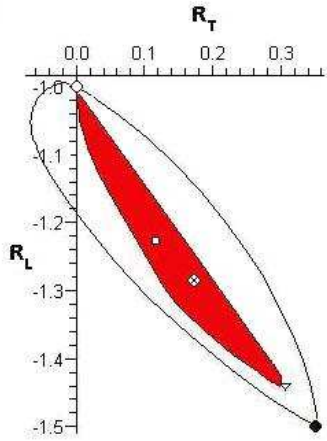}}
\hfill
\subfigure[$\lam = 1$]{\includegraphics[width=.40\textwidth]{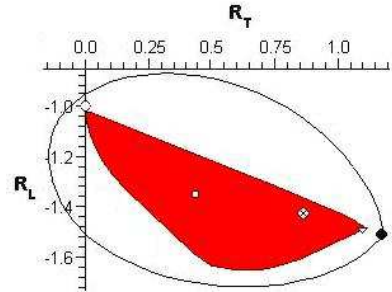}}
\end{center}
\caption{The sets $\mathcal{R}_\lam$ with $\lam = 0^+,\, 0.1,\, 0.3,\, 1$ are shown here separately.
The values $R[\nu, \lam]$, with $\nu = \nu_0,\, \nu_\star,\, \nu_{\text{rect}},\, \nu_\triangledown,\, \nu_\otimes$ are indicated by the signs $\circ$,
\mbox{\footnotesize$\bullet$}, \mbox{\tiny$\Box$}, \mbox{\footnotesize$\triangledown$}, \mbox{\footnotesize$\otimes$},
respectively.
On Fig.\,7(a), the region   generated by triangular cavities is marked by a (slightly shifted) vertical dashed line.
It is the interval with the endpoints $(0, -1)$ and $(0, -1.42)$. On Fig.\,7(b)--7(d),  the regions generated by triangular cavities are colored red.}
\end{figure}

\newpage

\begin{figure}[htb!]
\begin{center}
\subfigure[][]{\includegraphics[width=0.48\textwidth]{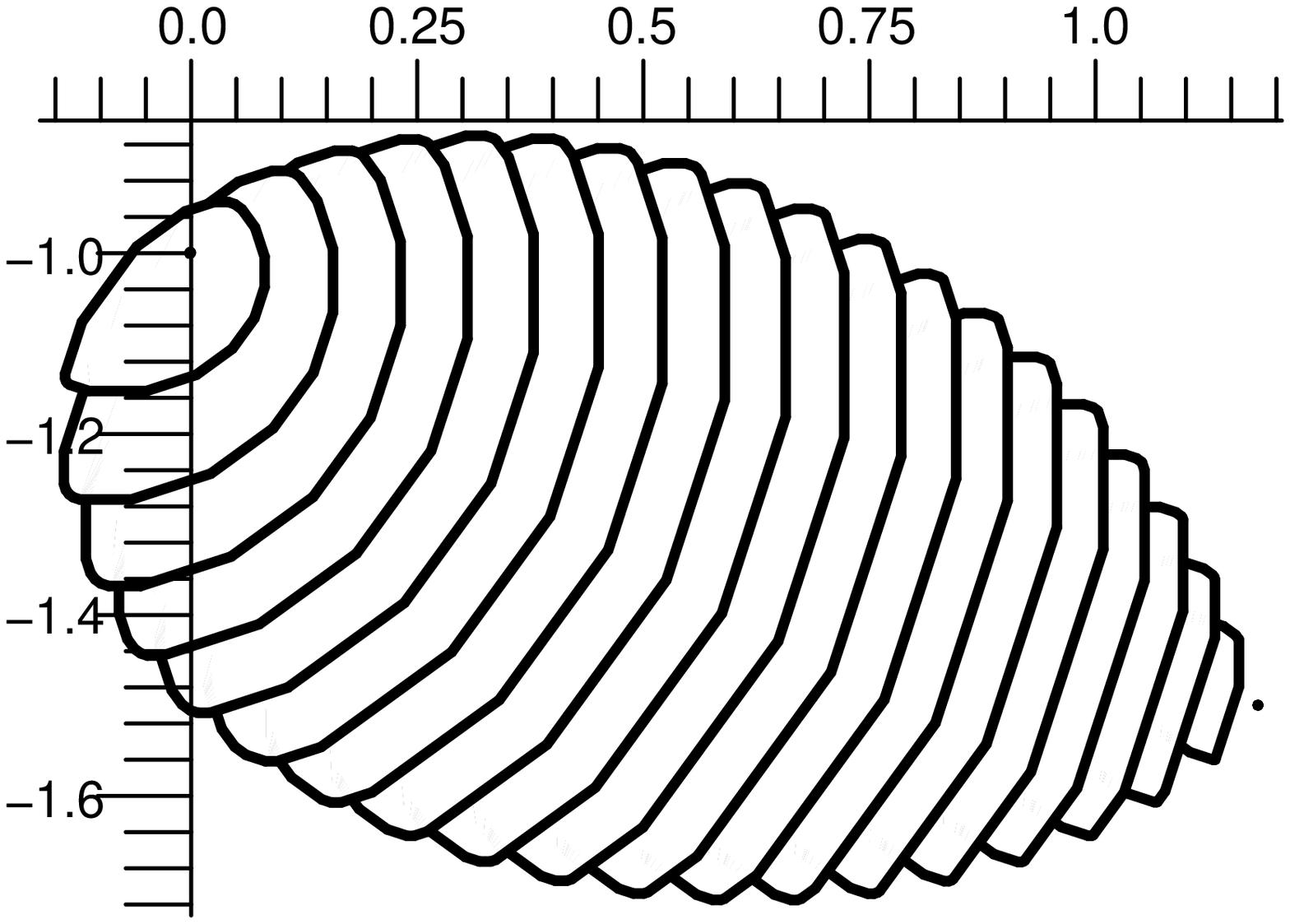}}
\hfill
\subfigure[][]{\includegraphics[width=0.48\textwidth]{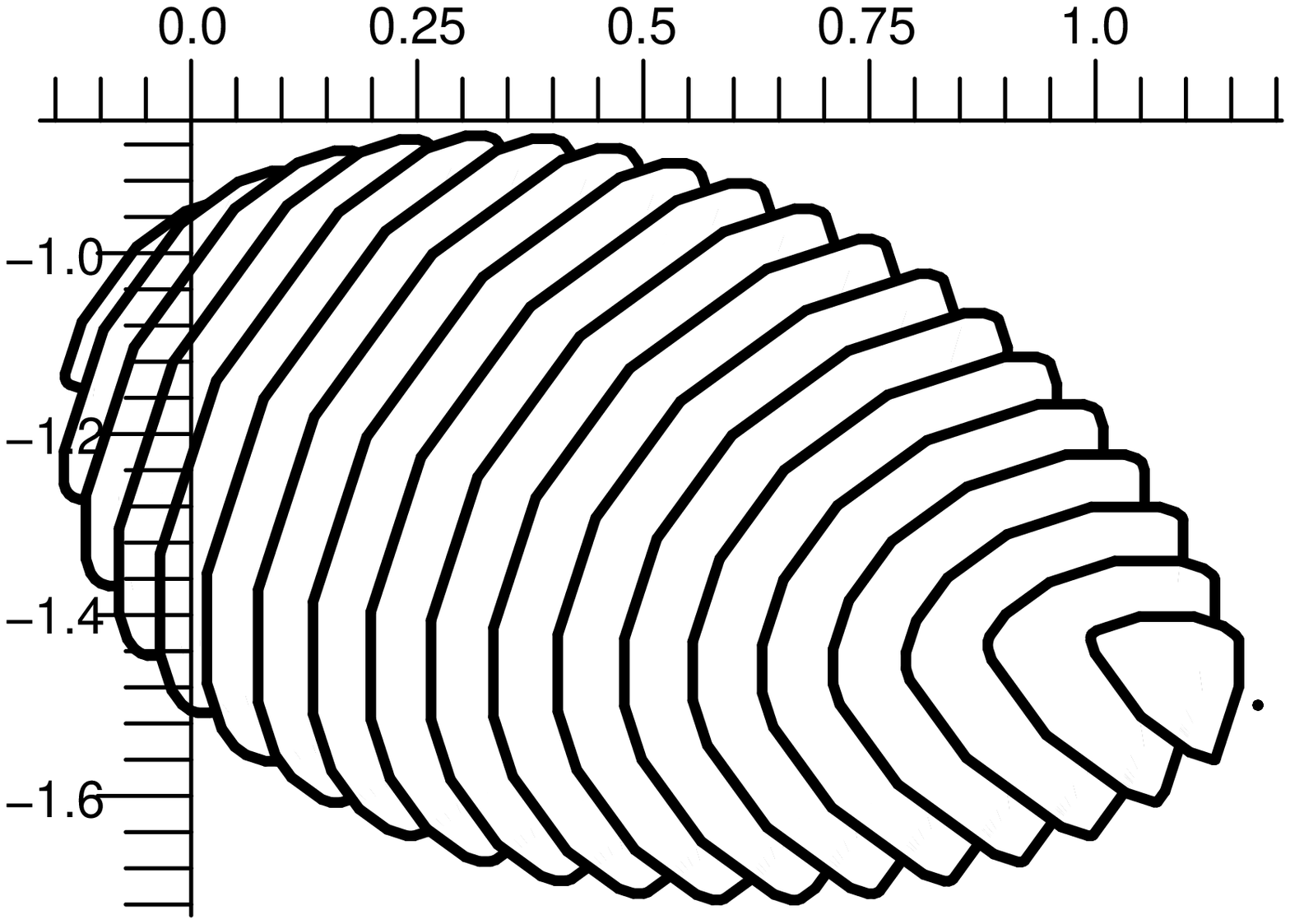}}
\end{center}
\caption{The ``level sets'' $\mathcal{R}_{1,c} = \{ \, \text{all possible values of } \vec R[\nu, 1], \text{ with } R_I[\nu, 1] = c \, \}$ are shown for 21 values of $c$, from left to right: $\,c = 0,\, -0.075,\, -0.15,$ $
-0.225, \ldots, -1.425$,\, $-1.5$.
Figures 8(a) and 8(b) represent, respectively, ``view from above'' and ``view from below'' on these sets.}
\end{figure}

\section{Dynamics of a rough disk}

The motion of a spinning rough {disk} $B$ is determined by the values $R_T[\nu_B,\lam]$,\, $R_L[\nu_B,\lam]$, and $R_I[\nu_B,\lam]$. For the sake of brevity, below we omit the fixed argument $\nu$ and write $R(\lam)$ instead of $R[\nu,\lam]$. Recall that the absolute value of the  disk velocity is denoted by $v = |\vec v|$, and the angular velocity equals $\om = \lam v/r$. Denote by $\theta$ the angle the velocity makes with a fixed direction in an inertial frame of reference.

Using (\ref{resistance}) and (\ref{moment of resistance}), one rewrites the equations of motion (\ref{dynam1}), (\ref{dynam2}) in the form
$$
\frac{dv}{dt}\ =\ \frac{8r\rho v^2}{3M}\, R_L(\lam),
$$
$$
\frac{d\theta}{dt}\ =\ -\frac{8r\rho v}{3M}\, R_T(\lam),
$$
$$
\frac{d(\lam v)}{dt}\ =\ \frac{8r^3\rho v^2}{3I}\, l R_I(\lam).
$$

Recall that $\relmom = Mr^2/I$  is the inverse relative moment of inertia. The special values that $\beta$ can take are: $\relmom = 1$, when the mass is concentrated on the disk boundary, and $\relmom = 2$, when
the mass is uniformly distributed inside the disk. In the intermediate case, when the mass is arbitrarily (generally speaking, non-uniformly) distributed inside the  disk, it holds $\relmom \ge 1$.

With the help of the change of variables
\begin{equation}\label{tau}d\tau = \frac{8r\rho v}{3M}\, dt\end{equation}
the equations above are transformed into the following relations:
 \beq\label{motion 2}
 \frac{d\lam}{d\tau}\ =\ \relmom R_I(\lam) - \lam R_L(\lam),
 \eeq
 \beq\label{motion 1}
 \frac{dv}{d\tau}\ =\ v\, R_L(\lam), \quad \quad
 \eeq
 \beq\label{motion 3}
 \frac{d\theta}{d\tau}\ =\ -R_T(\lam). \quad \quad
 \eeq

Denote by $s$ the path length of the disk; thus, $ds/dt = v$. From (\ref{tau}) one readily finds that $s$ is proportional to $\tau$,\, $s = \frac{3M}{8r\rho}\, \tau$.

Below we solve the system of equations (\ref{motion 2})--(\ref{motion 3}) for the cases 1 -- 3 considered in subsection 3.1. Next, we determine the dynamics numerically for some kinds of triangular cavity. A more detailed study of dynamics in the general case will be addressed elsewhere.

\begin{enumerate}

\item {\it Circle}.~ One has $d\lam/d\tau = -\lam$,\, $dv/d\tau = -v$ and $d\theta/dt = 0$;  therefore the circle moves
straightforward. Solving these equations, one concludes that its center moves according to  the equation $\vec x(t) = \frac{3M}{8r\rho}\, \ln(t - t_0) \vec e + \vec x_0$, where $t_0 \in \RRR$,\, $\vec e \in S^1$  and $\vec x_0 \in \RRR^2$ are constants. Thus, having started the motion at some moment, the circle passes a half-line during an infinite time. This equation also indicates that the motion cannot be defined for all values of $t \in \RRR$.

\item {\it Retroreflector}.~ Here the system  (\ref{motion 2})--(\ref{motion 3}) takes the form
  \beq\label{motion 4}
d\lam/d\tau = -3\lam(\relmom - 1)/2, \ \ dv/d\tau = -3v/2, \ \ d\theta/d\tau = -3\pi\lam/8.
 \eeq
In the case  $\relmom = 1$,  one evidently  has $\lam = $ const. The disk moves along a circumference of radius $M/(\pi r\rho\lam)$ in the direction opposite to the angular velocity of rotation: if the disk rotates counterclockwise then its center moves clockwise along the circumference. The radius of the circumference is proportional to the disk mass and inversely proportional to the relative angular velocity. The path length is proportional to the logarithm of time, $s(t) = \frac{M}{4r\rho}\, \ln(t - t_0)$.

In the case $\relmom > 1$, we have $s(t) = \frac{M}{4r\rho} \ln(t - t_0)$,\, $\theta = \theta_0 + \text{const} \cdot \exp(-(\relmom - 1)\frac{4r\rho}{M}\, s)$, and $\lam = \frac{4}{\pi} (\relmom - 1) (\theta - \theta_0)$.  The path length once again depends logarithmically on the time, the relative angular velocity $\lam$ converges to zero, and the direction $\theta$ converges to a limiting value $\theta_0$;  thus, the values $\lam$ and $\theta$ are exponentially decreasing functions of the path length and are inversely proportional to the $(\bt - 1)$th degree of the time passed since a fixed moment. The trajectory of motion is a semibounded curve that approaches an asymptote as $t \to +\infty$.

\item {\it Rectangular cavity}.~ Equations of motion  (\ref{motion 2})--(\ref{motion 3}) in this case take the form
$$
d\lam/d\tau = -3\lam(\relmom - 5/3)/4, \ \ dv/d\tau = -5v/4, \ \ d\theta/d\tau = -3\pi\lam/16.
$$
Solving these equations, one obtains
$\tau = \frac{4}{5} \ln(t - t_0)$,\,
$v = v_0\, e^{-5\tau/4}$,\, $\lam = \lam_0\, e^{3\tau(5/3 - \relmom)/4}$,\, $\theta = \theta_0
+ \frac{\pi \lam_0}{4(\relmom - 5/3)}\, e^{3\tau(5/3 - \relmom)/4}$. Thus, the {path} depends on $t$ logarithmically, and the relative angular velocity and the rotation angle are proportional to $(t-t_0)^{1-3\relmom/5}$ and to $\exp \big( {\frac{2r\rho(5-3\relmom)}{3M}\, s} \big)$.

If $\relmom < 5/3$ then $\lam$ and $\theta$ tend to infinity, and the trajectory of the disk center is a converging spiral. In the case $\relmom > 5/3$, $\lam$ converges to zero, $\theta$ converges to a constant value, and the trajectory is a semibounded curve approaching an asymptote as $t \to +\infty$. In the case $\relmom = 5/3$,\, $\lam$ is constant, and the trajectory is a circumference of radius $2M/(\pi r\rho\lam)$.

\end{enumerate}

Finally, let us describe numerically some triangular cavities. It is helpful to consider the function $g(\lam) = \lam R_L(\lam)/R_I(\lam)$ and rewrite the equation (\ref{motion 2}) in the form
 \beq\label{motion lam}
\frac{d\lam}{d\tau}\ =\ -R_I(\lam) (g(\lam) - \relmom).
 \eeq

The function $g(\lam)$ has been calculated for several kinds of cavity. On Fig.\,9\,(a), the two cases are illustrated where the cavity is an isosceles triangle with the angles (i) $30^0$,\, $120^0$,\, $30^0$ and (ii) $60^0$,\, $60^0$,\, $60^0$. We see that $g(\lam)$ monotonically increases in the case (i) and has three intervals of monotonicity in the case (ii). In both cases, it is valid $R_I(\lam) < 0$. This implies, in the case (i), that the disk trajectory is a converging spiral, if $\relmom < 1.5$. If $\relmom > 1.5$, the trajectory depends on the initial conditions and may take the form of a converging spiral or a curve approaching a straight line.

\vspace{7mm}
\begin{figure}[htb!]
\begin{center}
\subfigure[][]{\includegraphics[width=.48\textwidth]{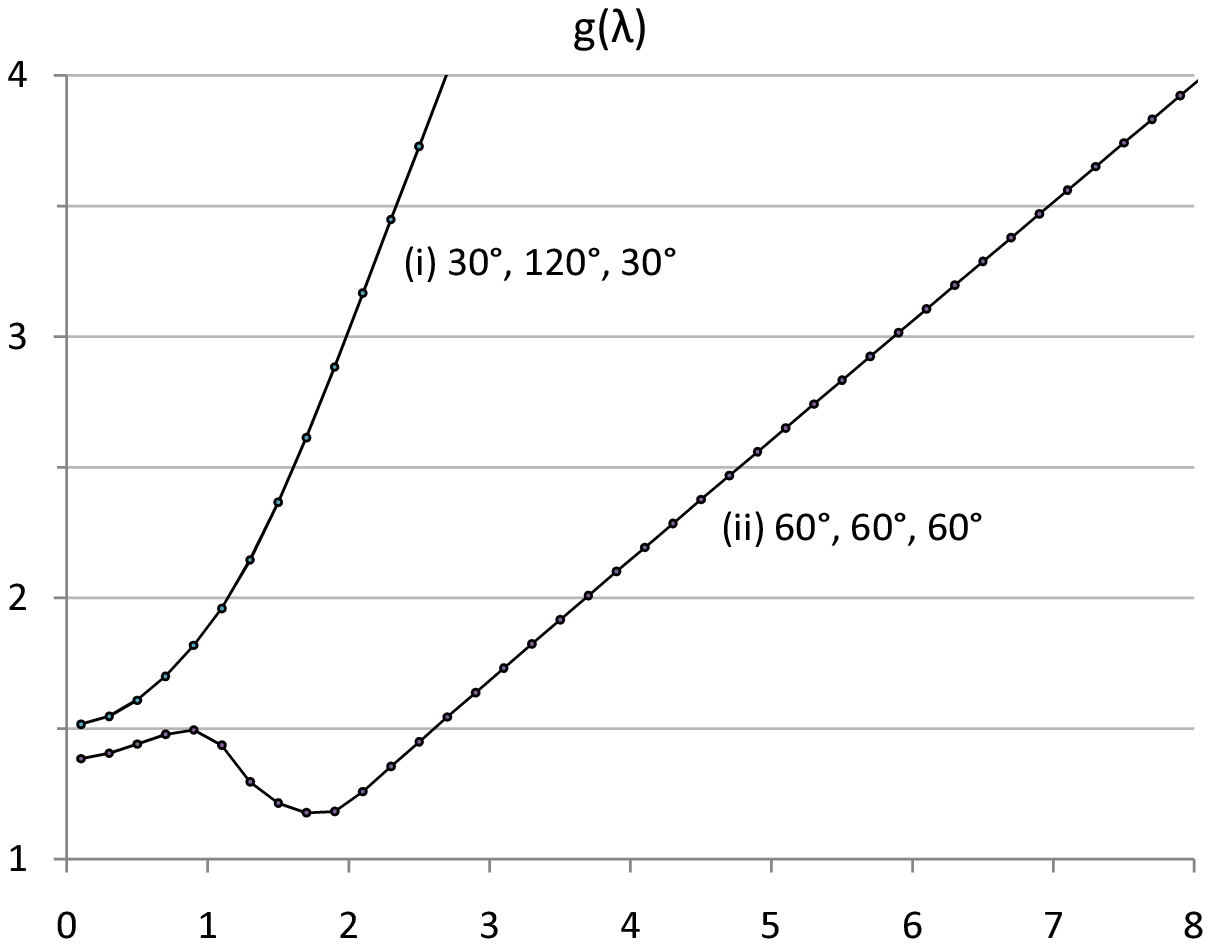}}
\hfill
\subfigure[][]{\includegraphics[width=.48\textwidth]{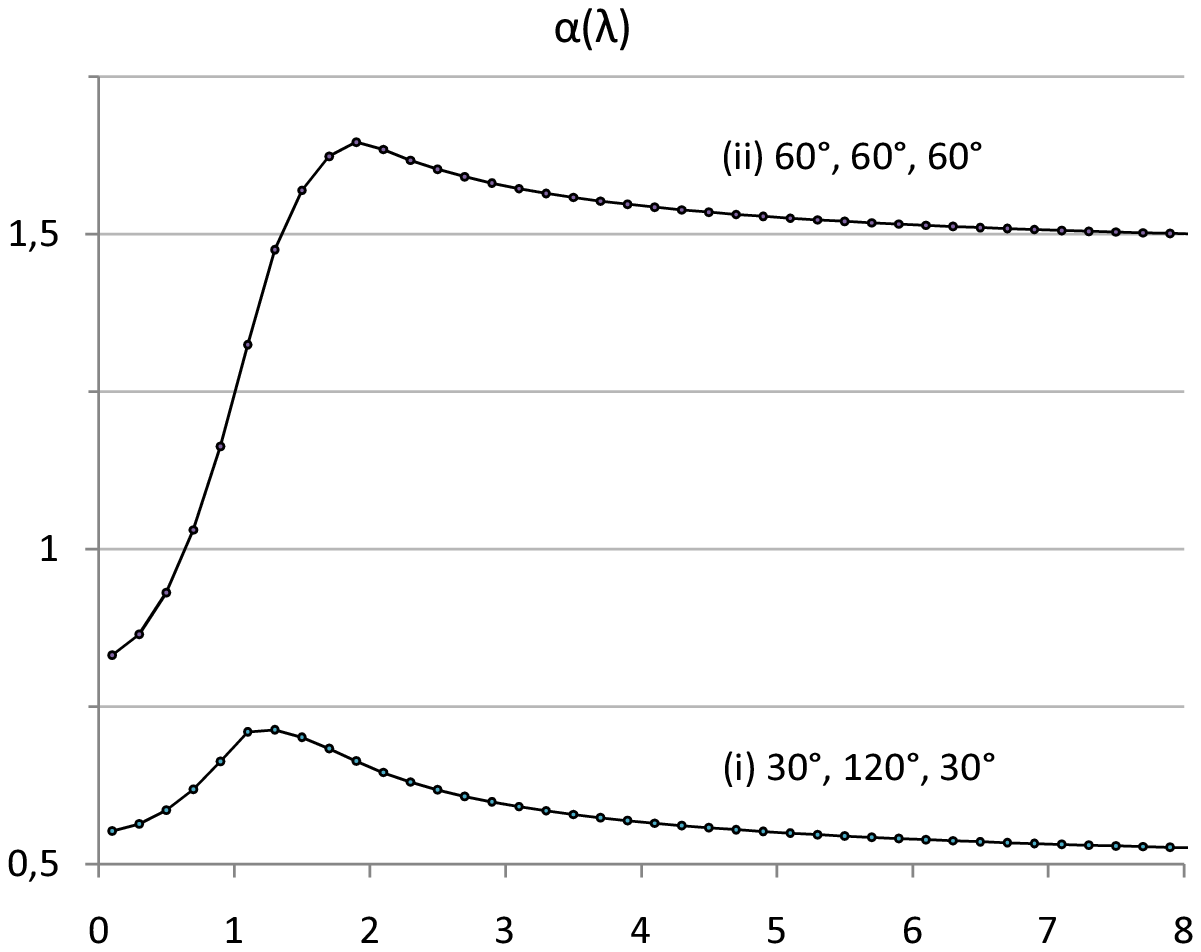}}
\end{center}
\caption{ The functions (a) $g(\lam) =  \lam R_L(\lam)/R_I(\lam)$ and (b) $\al(\lam) = 4R_T(\lam)/\lam$ are shown for the triangular cavities with the angles (i) $30^0$,\, $120^0$,\, $30^0$ and (ii) $60^0$,\, $60^0$,\, $60^0$.}
\end{figure}
\vspace{7mm}

The disk behavior is  more rich in the case (ii) of equilateral triangle. If $1.38 < \relmom < 1.49$ then three kinds of asymptotic behavior may be realized, depending on the initial conditions: (I) the trajectory is a converging spiral; (II) the trajectory approaches a circumference; and (III) the trajectory approaches a straight line  (see Fig.\,10). If $1.16 < \relmom < 1.38$, only two asymptotic behaviors of the types (I) and (II) are possible; if $\relmom > 1.49$ then the possible behaviors are (I) and (III); and finally, if $\relmom < 1.16$, the asymptotic behavior is always (I).
 \newpage

\begin{figure}[!h]
\begin{center}
\includegraphics[width=.9\textwidth]{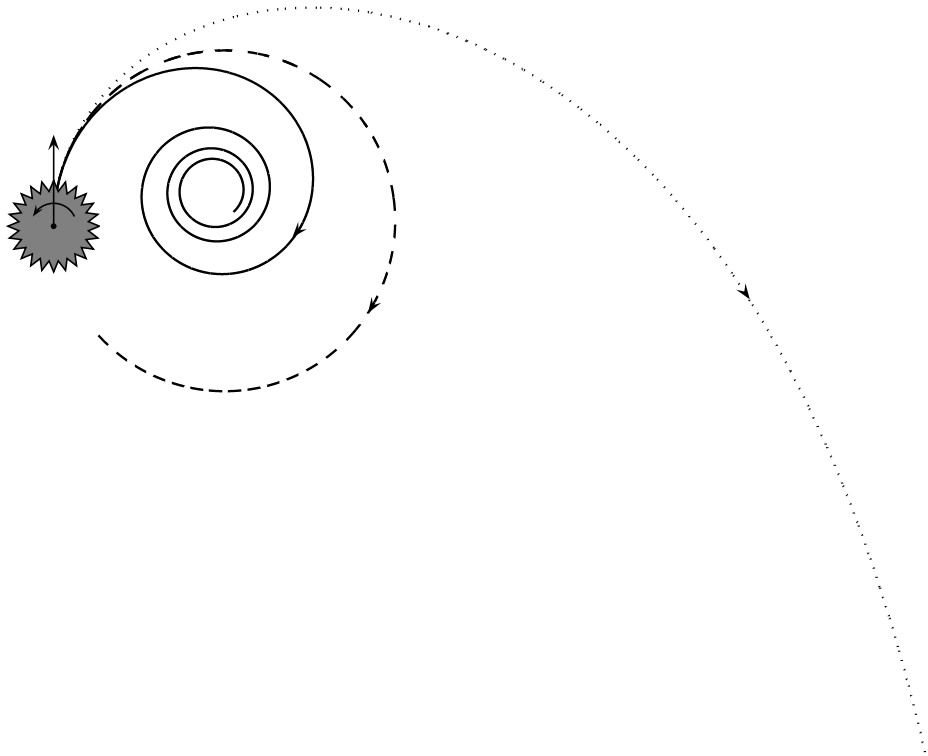}
\end{center}
\caption{Three kinds of asymptotic behavior of a rough disk with roughness formed by equilateral triangles  and with $1.38 < \relmom < 1.49$: (I) converging spiral (solid line); (II) circumference (dashed line); (III) curve approaching a straight line (dotted line).}
\end{figure}

In the case of triangular cavities, as our numerical evidence shows, the function $g(\lam)$ monotonically increases  for $\lam$ sufficiently large and $\lim_{\lam\to+\infty} g(\lam) = +\infty$. This implies that the trajectory is a converging spiral for appropriate initial conditions  (namely, if the initial angular velocity is large enough). If, besides, $\relmom$ is large enough (that is, the mass of the disk is concentrated near the center), the trajectory may also be a curve approaching a straight line. If the function $g$ admits intervals of monotone decrease (as for the case of equilateral triangle) then the trajectory may approach a circumference. The length of the disk path is always proportional to the logarithm of time.

\section{Conclusions and comparison with the previous works}

In the papers \cite{BSE}, \cite{IY}, \cite{W}, and \cite{WH}, the same problem as in the present paper is studied: determine the force acting on a spinning body moving through a rarefied gas. In \cite{IY}, additionally, the moment of this force slowing down the body's rotation is calculated.  The following shapes of the body have been considered: a sphere, a cylinder \cite{IY,WH}, convex bodies of revolution \cite{IY}, and right parallelepipeds of regular polygon section \cite{WH}. The interaction of the body with the gas particles is in general nonelastic:  a fraction $1 - \al_\tau$ of the incident particles is elastically reflected according to the rule ``the angle of incidence is equal to the angle of reflection'', while the fraction $\al_\tau$ of the particles reaches thermal equilibrium with the body's surface, and is reflected as a Maxwellian \cite{BSE},\cite{IY},\cite{W}. In the paper \cite{WH}, a somewhat different model of interaction is considered, where  the reflected particles acquire a fraction $\al_\tau$ of tangential momentum of the rotating body. The transversal force results from the tangential friction and acts on the body in the direction corresponding to the reverse Magnus effect. It is remarkable that for different models and for different shapes of the body, the formula for the transversal force is basically the same. If the rotation axis is perpendicular to the direction of the body's motion then this force equals $\frac 12\, \al_\tau \massgas \om v$, where $\massgas$ is the mass of the gas displaced by the body, $\om$ is the angular velocity of the body, and $v$ is its translation velocity. (Note that in \cite{BSE} this formula appears in the limit of infinite heat conductivity or zero gas temperature.) In \cite{WH}  it was found that for parallelepipeds of regular $n$-gon section with  $n$ odd, transversal force depends on time, and the value of this force was determined. It is easy to calculate, however, that the time average of the force is equal to the same quantity $\frac 12\, \al_\tau \massgas \om v$.

In the present paper, to the contrary, it is assumed that all collisions of particles with the body are absolutely elastic (that is, $\al_\tau = 0$), and therefore there is no tangential friction. We restrict ourselves to the two-dimensional case and suppose that the body is a disk with small cavities on its boundary, or a {\it rough} disk. The gas has zero temperature. The Magnus effect is due to multiple reflections of particles in the cavities.
According to (\ref{resistance}), the transversal force equals $\frac 12\, \al(\lam)\, \massgas \om v$, where $\lam = \om r/v$,\, $\massgas = \pi r^2 \rho$,\, $\al(\lam) = \al(\lam, \nu) = (16/3\pi) R_T[\nu, \lam]/\lam$, and $\nu$ is the measure characterizing the shape of the cavities.  Thus, here $\massgas$ also means the total mass of gas particles displaced by the body. The function $\al$ depends on both $\nu$ and $\lam$. In particular, $\al$ varies between $-0.409$ and 2 for $\lam = 0.1$, between $-0.378$ and 2 for $\lam = 0.3$, and between $-0.248$ and 2 for $\lam = 1$. We conjecture that $\lim_{\lam\to\infty} \inf_\nu \al(\lam,\nu) = 0$ and $\lim_{\lam\to\infty} \sup_\nu \al(\lam,\nu) = 2$. The graphs of the function $\al(\lam)$ with $\nu$ corresponding to triangular cavities with the angles (i) $30^0$,\, $120^0$,\, $30^0$ and (ii) $60^0$,\, $60^0$,\, $60^0$ are shown on Fig.\,9(b). We see that this function significantly depends on the velocity of rotation $\lam$; in general, the variation of $\al(\lam)$ with $\nu$ fixed can be more that twofold.

 \vspace{3mm}

In our opinion, the reverse Magnus effect appears due to the following two factors: (i) nonelastic interaction of particles with the body, and (ii) multiple reflections from cavities on the body's surface. In the papers \cite{BSE},\cite{IY},\cite{W},\cite{WH} the influence of the factor (i) has been studied, while we concentrate on the factor (ii). Our conclusion is that the influence of both factors is unidirectional, so they strengthen each other. Moreover, the formulas for the transversal force are similar; one should just substitute the function $\al(\lam, \nu)$ for $\al_\tau$. We have seen that $\al(\lam, \nu)$ can be significantly greater than 1, while $\al_\tau \le 1$. This can be just an artefact of two-dimensionality (in our case). In a forthcoming work we are planning to extend our consideration to the three-dimensional case and to consider media with non-zero temperature.

\section*{Acknowledgements}

This work was supported by {\it Centre for Research on Optimization and Control} (CEOC) from the ``{\it Fundação para a Ciência e a Tecnologia}'' (FCT), cofinanced by the European Community Fund FEDER/POCTI, and by FCT: research project PTDC/MAT/72840/2006.

The authors are grateful to John Gough for useful comments on the text and to
Gennady Mishuris for the help in drawing Figure 10.


\begin{thebibliography}{99}

\bibitem{BS}
 K.\,I. Borg and L.\,H. Söderholm. {\it Orbital effects of the Magnus force on a
spinning spherical satellite in a rarefied atmosphere}. Eur. J. Mech. B/Fluids {\bf
27}, 623-631 (2008).

\bibitem{BSE}
K.\,I. Borg, L.\,H. Söderholm and H. Essén. {\it Force on a
spinning sphere moving in a rarefied gas}. Physics of Fluids {\bf
15}, 736-741 (2003).

\bibitem{IY}
S.\,G. Ivanov and A.\,M. Yanshin. {\it Forces and moments acting on
bodies rotating around a symmetry axis in a free molecular flow}.
Fluid Dyn. {\bf 15}, 449 (1980).

\bibitem{Mehta}
R.\,D. Mehta. {\it Aerodynamics of sport balls}. Annu. Rev. Fluid Mech. {\bf 17}, 151-189 (1985).

\bibitem{MatSb 04: mean resist}
A. Plakhov. {\it Newton's problem of the body of minimum mean
resistance}. Sbornik: Math. {\bf 195}, 1017-1037 (2004).

\bibitem{Nonlinearity}
A. Plakhov and P. Gouveia. {\it Problems of maximal mean resistance
on the plane}. Nonlinearity {\bf 20}, 2271-2287 (2007).

\bibitem{Arch Rat Mech}
A. Plakhov. {\it Billiards and two-dimensional problems of optimal resistance}.
Arch. Ration. Mech. Anal., DOI 10.1007/s00205-008-0137-1, 33 pp.

\bibitem{rough body 2D}
A. Plakhov. {\it Billiard scattering on rough sets: Two-dimensional case}.
SIAM J. Math. Anal.{\bf 40}, 2155-2178 (2009). DOI. 10.1137/070709700.

\bibitem{Prandtl}
 L. Prandtl. {\it Application of the ``Magnus Effect'' to the wind propulsion of ships}.
Die. Naturwissenschaft {\bf 13}, 93-108; transl. NACA-TM-367, June 1926.



\bibitem{RK}
 S.\,I. Rubinov and J.\,B. Keller. {\it The transverse force on a spinning sphere moving in a viscous fluid}. J. Fluid Mech., {\bf 11}, 447-459 (1961).

\bibitem{W}
C.-T. Wang. {\it Free molecular flow over a rotating sphere}.
AIAA J. {\bf 10}, 713 (1972).

\bibitem{WH}
P.\,D. Weidman and A. Herczynski. {\it On the inverse Magnus effect
in free molecular flow}. Physics of Fluids {\bf 16}, L9-L12 (2004).












\end{thebibliography}
\end{document}